\documentclass[12pt,reqno]{amsart}
\usepackage{a4,amsmath,mathrsfs}

\makeatletter \@addtoreset{equation}{section}

\makeatother \makeatletter
\newcommand{\eps}{\varepsilon}

\newcommand{\rn}[1]{{\bf R}^{#1}}
\newcommand\R{{\rn{}}}
\newcommand\Z{{\bf Z}}
\newcommand\N{{\bf N}}

\newcommand{\segop}{[\kern-2pt[}
\newcommand{\segcl}{]\kern-2pt]}
\newcommand{\Leb}[1]{{\mathscr L}^{#1}}
\newcommand{\Haus}[1]{{\mathscr H}^{#1}}
\newcommand{\res}{\mathop{\hbox{\vrule height 7pt width .5pt depth 0pt
\vrule height .5pt width 6pt depth 0pt}}\nolimits}

\newcommand{\modp}{{\rm mod}(p)}

\newcommand{\fflat}{{\mathscr F}}

\newcommand{\rc}[2]{{\mathcal I}_{{#1}}({#2})}      

\newcommand{\fc}[2]{{\mathscr F}_{{#1}}({#2})}       
\newcommand{\fcp}[2]{{\mathscr F}_{p,{#1}}({{#2}})}

\newcommand{\ic}[2]{{\bf I}_{{#1}}({#2})}       

\newcommand\bdry{\partial}
\newcommand\mass{{\bf M}}
\newcommand\Lbrack{\segop}
\newcommand\Rbrack{\segcl}
\newcommand{\on}[1]{|_{#1}}

\newtheorem{definition}{Definition}[section]
\newtheorem{lemma}[definition]{Lemma}
\newtheorem{theorem}[definition]{Theorem}
\newtheorem{proposition}[definition]{Proposition}
\newtheorem{corollary}[definition]{Corollary}
\newtheorem{remark}[definition]{Remark}

\title[Rectifiability of flat $\Z_p$ chains]
{Rectifiability of flat chains in Banach spaces with coefficients in
$\Z_p$}

\author{ Luigi Ambrosio}

\address{Scuola Normale Superiore, Piazza Cavalieri 7, 56100 Pisa, Italy}
\email{l.ambrosio@sns.it}

\author{ Stefan Wenger}
\thanks{The first author's work was partially supported by a MIUR PRIN06 grant. The second author's work was partially supported by NSF grant DMS 0707009.}
\address{University of Illinois at Chicago, 851 S. Morgan Street, Chicago, IL 60607-7045, USA}

\email{wenger@math.uic.edu}

\begin{document}

\maketitle

\section{Introduction}

Aim of this paper is a finer analysis of the group of flat chains
with coefficients in $\Z_p$ introduced in \cite{ambkatz}, by taking
quotients of the group of integer rectifiable currents, along the
lines of \cite{Z,federer}. We investigate the typical questions of
the theory of currents, namely rectifiability of the
measure-theoretic support and boundary rectifiability. Our main
result can also be interpreted as a closure theorem for the class
of integer rectifiable currents with respect to a (much) weaker
convergence, induced by flat distance $\modp$, and with respect to
weaker mass bounds. A crucial tool in many proofs is the
isoperimetric inequality proved in \cite{ambkatz} with universal
constants.

In order to illustrate our results we start with a few basic
definitions. Let us denote by $\rc{k}{E}$ the class of integer
rectifiable currents with finite mass in a metric space $E$ and let
us given for granted the concepts of boundary $\partial$, mass ${\bf
M}$, push-forward in the more general context of currents (see
\cite{ak2} and the short appendix of \cite{ambkatz}). We denote by
$\fc{k}{E}$ the currents that can by written as $R+\partial S$ with
$R\in\rc{k}{E}$ and $S\in\rc{k+1}{E}$. It is obviously an additive
Abelian group and
\begin{equation}\label{amb2}
T\in\fc{k}{E}\qquad\Longrightarrow\qquad\partial T\in\fc{k-1}{E}.
\end{equation}
$\fc{k}{E}$ is a metric space when endowed with the the distance
$d(T_1,T_2)=\fflat(T_1-T_2)$, where
$$
\fflat(T):=\inf\left\{{\bf M}(R)+{\bf M}(S):\ R\in\rc{k}{E},\,\,
S\in\rc{k+1}{E},\,\,T=R+\partial S\right\}.
$$
The subadditivity of $\fflat$, namely $\fflat(nT)\leq n\fflat(T)$,
ensures that $d$ is a distance, and the completeness of the groups
$\rc{k}{E}$, when endowed with the mass norm, ensures that
$\fc{k}{E}$ is complete. Also, the boundary rectifiability theorem
in $\rc{k}{E}$ yields
\begin{equation}\label{comedia}
\left\{T\in\fc{k}{E}:\ {\bf M}(T)<\infty\right\}=\rc{k}{E}.
\end{equation}
For $T\in\fc{k}{E}$ we define:
$$
\fflat_p(T):=\inf\left\{\fflat(T-pQ):\ Q\in\fc{k}{E}\right\}.
$$
The definition of $\fflat$ gives
$$
\fflat_p(T)=\inf\left\{{\bf M}(R)+{\bf M}(S):\ T=R+\partial
S+pQ,\,\,R\in\rc{k}{E},\,\,
S\in\rc{k+1}{E},\,\,Q\in\fc{k}{E}\right\}.
$$
Obviously $\fflat_p(T)\leq\fflat(T)$ and therefore we can introduce
an equivalence relation $\modp$ in $\fc{k}{E}$, compatible with the
group structure, by saying that $T=\tilde{T}$ $\modp$ if
$\fflat_p(T-\tilde{T})=0$. Our main object of investigation will be
the spaces
$$
\fcp{k}{E}:=\left\{[T]:\ T\in\fc{k}{E}\right\}.
$$
The equivalence classes are closed in $\fc{k}{E}$ and (see
\eqref{amb3} in the next section) the boundary operator can be
defined also in the quotient spaces $\fcp{k}{E}$ in such a way that
$$
\partial [T]=[\partial T]\in\fcp{k-1}{E}\qquad\forall T\in\fc{k}{E}.
$$

In $\fc{k}{E}$ one can also define a (relaxed) notion of $p$-mass
${\bf M}_p$ by
\begin{equation}\label{mp}
{\bf M}_p(T):=\inf\left\{\liminf_{h\to\infty}{\bf M}(T_h):
T_h\in\fc{k}{E},\,\,\fflat_p(T_h-T)\to 0\right\}.
\end{equation}
Since ${\bf M}_p(T)={\bf M}_p(T')$ if $T=T'$ $\modp$ the definition
obviously extends to the quotient spaces $\fcp{k}{E}$. As in the
standard theory of currents, a local variant of this definition
provides a $\sigma$-additive Borel measure, that we shall denote by
$\|T\|_p$, whose total mass is $\mass_p(T)$.

>From now on, we shall assume that $E$ is a compact convex subspace
of a Banach space $F$ and a Lipschitz retract of it. In addition, we shall assume that $F$ satisfies a strong finite-dimensional
approximation property (precisely stated in Definition~\ref{dsfda})
that covers, for instance, all Hilbert spaces.

We can now state the main result of this paper.

\begin{theorem}[Rectifiability of flat chains $\modp$]\label{tmain1}
If $T\in\fc{k}{E}$ has finite ${\bf M}_p$ mass, then $\|T\|_p$ is
concentrated on a countably $\Haus{k}$-rectifiable set with finite
$\Haus{k}$-measure.
\end{theorem}

We don't know whether the result is true without the
finite-dimensional approximation assumption, unless $k=0,1$. In
general, without this assumption, we are able to prove
rectifiability only of a the ``slice mass'' $\|T\|_p^*$ (see
Definition~\ref{dslicemass} and \eqref{realfl1}) built using the
$0$-dimensional slices of the flat chain, and the validity in
general spaces of the equality $\|T\|_p=\|T\|_p^*$ is still an open
problem.

Since $\partial$ maps $\fc{k}{E}$ into $\fc{k-1}{E}$, the next
result is a direct consequence of Theorem~\ref{tmain1}.

\begin{corollary}[Boundary rectifiability]\label{cmain2}
If $T\in\fc{k}{E}$ and if $\partial T$ has finite ${\bf M}_p$-mass,
then $\|\partial T\|_p$ is concentrated on a countably
$\Haus{k-1}$-rectifiable set with finite $\Haus{k-1}$-measure.
\end{corollary}

Notice that in Corollary~\ref{cmain2} finiteness of mass of $T$ is
not needed. As a corollary we obtain an extension $\modp$ of
\eqref{comedia}, namely flat chains $\modp$ with finite ${\bf M}_p$
mass coincide with equivalence classes of integer rectifiable
currents. These classes have been considered in \cite{ambkatz} in
connection with isoperimetric and filling radius inequalities.

\begin{corollary}\label{cmain3}
If $T\in\fc{k}{E}$ has finite ${\bf M}_p$ mass, then there exists
$S\in\rc{k}{E}$ with $S=T$ $\modp$. In addition, $S$ can be chosen
so that $\mass_p(T)=\mass(S)$.
\end{corollary}

We give a detailed proof of Corollary~\ref{cmain3} at the end of the
paper. We obtain also, as a byproduct, the following closure theorem
for $\rc{k}{E}$: in comparison with the results in \cite{ak2} the
$\fflat_p$ convergence (instead of the weak convergence in the
duality with all Lipschitz differential forms), and the bounds only
on the ${\bf M}_p$ mass (instead of the stronger mass bounds)
are considered. Obviously the result can also be stated
as a closure theorem in $\fcp{k}{E}$.

\begin{corollary}[Closure theorem]\label{cmain4}
Assume that $(T_n)\subset\rc{k}{E}$ satisfies $\fflat_p(T_n-T)\to 0$
for some $T\in\fc{k}{E}$. If $\sup_n{\bf M}_p(T_n)<\infty$,
then there exists $S\in\rc{k}{E}$ with $S=T$ $\modp$.
\end{corollary}

We conclude the introduction with a short plan of the paper. In
Section~\ref{s1} we recall the basic results we need on flat chains
and flat chains $\modp$, borrowing some results from \cite{ambkatz}.
In Section~\ref{s2} we study more in detail the slicing operator and
the measure $\|T\|_p$. The main result is that a flat chain with
finite mass and boundary with finite mass is uniquely determined by
its slices. In this section we don't rely, as in \cite{white1} on
the use of the deformation theorem of \cite{white2}, not available
in our context. We heavily use, instead, the isoperimetric
inequality: in turn, this inequality (derived as well in
\cite{white1} as a consequence of the deformation theorem) is proved
in \cite{ambkatz} without using the deformation theorem. In
Section~\ref{s3} we make a finer analysis of 1-dimensional flat
chains $\modp$ and we provide a direct proof of their
rectifiability; this is a crucial ingredient to estabilish the
rectifiability of the slice mass $\|T\|_p^*$ of higher dimensional
flat chains, following basically the procedure in \cite{white1}.
This procedure is implemented in Section~\ref{s4} and
Section~\ref{s5} and leads to the proof that $\|T\|_p^*$ is
concentrated on a countably $\Haus{k}$-rectifiable set: the main
difference with respect to \cite{white1} consists in the fact that
the whole family of $1$-Lipschitz projections, instead of the
projections on the coordinate planes typical of the Euclidean case,
has to be considered. In this respect, notice that still a $BV$
estimate analogous to the one in \cite{ak1} is available in this
setting, see Remark~\ref{rBV}, and it is likely that also some
adaptations of the ideas in \cite{ak1} might provide a different
proof of the rectifiability of $\|T\|_p^*$. Finally, in
Section~\ref{s6} we complete our analysis getting a concentration
set with finite $\Haus{k}$-measure and proving the equality
$\|T\|_p=\|T\|_p^*$ in the class of spaces having the
finite-dimensional approximation property.

\section{Notation and basic results on flat chains}\label{s1}

We use the standard notation $B_r(x)$ for the open balls in $E$,
${\rm Lip}(E)$ for the space of Lipschitz real-valued functions
and ${\rm Lip}_b(E)$ for bounded Lipschitz functions.
Now we recall a few basic facts on flat chains and flat
chains $\modp$ mostly estabilished in \cite{ambkatz}.

Throughout the paper we assume that $E$ is a compact convex subset
of a Banach space. Denoting by $\ic{k}{E}$ the space
$$
\ic{k}{E}:=\left\{T\in\rc{k}{E}:\ \partial T\in\rc{k-1}{E}\right\},
$$
this assumption ensures the density of $\ic{k}{E}$ in $\fc{k}{E}$
(see \cite{ambkatz}), and this gives the possibility to extend the
restriction and slicing operators from $\ic{k}{E}$ to $\fc{k}{E}$.

\subsection{Inequalities}
Notice that
\begin{equation}\label{amb1}
\fflat(\partial T)\leq\fflat(T),\qquad\qquad \forall T\in\fc{k}{E}.
\end{equation}
In addition, since $\partial (\varphi_\sharp
S)=\varphi_\sharp(\partial S)$, the inequality $\mass(\varphi_\sharp
R)\leq [{\rm Lip}(\varphi)]^k\mass(R)$ for $R$ $k$-dimensional gives
\begin{equation}\label{amb3bis}
\fflat(\varphi_\sharp T)\leq [{\rm Lip}(\varphi)]^k\fflat(T)
\quad\text{for all $T\in\fc{k}{E}$, $\varphi\in {\rm Lip}(E,\R^k)$.}
\end{equation}
In addition, \eqref{amb1} together with \eqref{amb2} give
\begin{equation}\label{amb3}
\fflat_p(\partial T)\leq\fflat_p(T),\qquad\quad \forall
T\in\fc{k}{E},
\end{equation}
while \eqref{amb3bis} gives
\begin{equation}\label{amb3ter}
\fflat_p(\varphi_\sharp T)\leq [{\rm Lip}(\varphi)]^k\fflat_p(T)
\end{equation}
for all $T\in\fc{k}{E}$, $\varphi\in {\rm Lip}(E,\R^k)$. In
particular, the push-forward operator can be defined in the quotient
spaces in such a way as to commute with the equivalence relation
$\modp$. Using \eqref{amb3} and the inequalities
$\fflat_p\leq\mass_p\leq\mass$
it is also easy to check that
\begin{equation}\label{may2}
\fflat_p(T)=\inf\left\{\mass_p(R)+\mass_p(S):\ R\in\rc{k}{E},\,\,
S\in\rc{k+1}{E}\right\}.
\end{equation}

\subsection{The restriction operator}

Let $u\in{\rm Lip}(E)$. In \cite{ambkatz} it is proved that the limit
\begin{equation}\label{fasto}
\lim_{h\to\infty}T_h\res\{u<r\}
\end{equation}
exists in $\fc{k}{E}$ for $\Leb{1}$-a.e. $r\in\R$ whenever $T_h$
have finite mass and $\sum_h\fflat(T_h-T)<\infty$. By construction
the operator $T\mapsto T\res\{u<r\}$ is additive and this definition
is independent, up to Lebesgue negligible sets, on the chosen
approximating sequence $(T_h)$, provided the ``fast convergence''
condition $\sum_h\fflat(T_h-T)<\infty$ holds. The construction
provides also the inequality
\begin{equation}\label{cachan}
\int_m^{*\ell}\fflat(T\res\{u<r\})\,dr\leq (\ell-m+{\rm
Lip}(u))\fflat(T)\qquad\forall m,\,\ell\in\R,\,\,m\leq\ell,
\end{equation}
where $\int^*$ denotes the outer integral. It follows immediately
from the additivity of the restriction operator that
\begin{equation}\label{cachanbis}
\int_m^{*\ell}\fflat_p(T\res\{u<r\})\,dr\leq (\ell-m+{\rm
Lip}(u))\fflat_p(T)\qquad\forall m,\,\ell\in\R,\,\,m\leq\ell,
\end{equation}
so that the restriction operator can be defined in the quotient
spaces $\fcp{k}{E}$ in such a way that
\begin{equation}\label{cachanter}
[T]\res\{u<r\}=[T\res\{u<r\}]\qquad\text{for $\Leb{1}$-a.e.
$r\in\R$.}
\end{equation}

\subsection{${\bf M}_p$-mass and $\|T\|_p$-measure}

Recall that the mass measure $\|T\|$ of $T\in\rc{k}{E}$ is the
finite nonnegative Borel measure characterized by
$\|T\|(\{u<r\})=\mass (T\res\{u<r\})$ for all $u\in {\rm Lip}(E)$
and $r\in\R$. In \cite[Theorem~7.1]{ambkatz} the authors proved the
existence, for all $T\in\fc{k}{E}$ of finite ${\bf M}_p$-mass, of a
finite nonnegative Borel measure $\|T\|_p$ satisfying
$$
{\bf M}_p(T\res\{u<r\})=\|T\|_p(\{u<r\}) \qquad\text{for
$\Leb{1}$-a.e. $r\in\R$}
$$
for all $u\in {\rm Lip}(E)$. In addition, since $\|T\|_p$ arises in
the proof of that result as the weak limit of $\|T_n\|$, where
$T_n\in\fc{k}{E}$ are such that $\fflat_p(T_n-T)\to 0$ and ${\bf
M}(T_n)\to{\bf M}_p(T)$, we can pass to the limit as $n\to\infty$ in
the inequalities
$$
\fflat(T_n\res\{u<s\}- T_n\res\{u<r\})\leq\|T_n\|(\{r\leq u< s\})
 \qquad r<s,
$$
taking into account that \eqref{cachanbis} gives
$\fflat_p(T_n\res\{u<r\}-T\res\{u<r\})\to 0$ for $\Leb{1}$-a.e.
$r\in\R$, to obtain
\begin{equation}\label{fasta}
\fflat_p(T\res\{u<s\}- T\res\{u<r\})\leq {\rm Lip}(u)\|T\|_p(\{r\leq
u\leq s\})
 \qquad\forall r,\,s\in\R\setminus N,\,\,r\leq s
\end{equation}
with $N$ Lebesgue negligible (possibly dependent on $T$ and $u$).

Using this fact, for chains $T$ with finite ${\bf M}_p$-mass we can
give a meaning to the restriction $[T]\res C$, for all fixed closed
set $C\subset E$, as follows: let $\pi$ be the distance function
from $C$, and let $N$ be as in \eqref{fasta} with $u=\pi$. If
$r_i\downarrow 0$ and $r_i\notin N$ then $T\res\{\pi<r_i\}$ is a
Cauchy sequence with respect to $\fflat_p$ and we denote by $[T]\res
C\in\fcp{k}{E}$ its limit. The lower semicontinuity of ${\bf M}_p$
provides also the inequality
$$
{\bf M}_p([T]\res C)\leq \|T\|_p(C).
$$
An analogous procedure (considering the sets $\{d(\cdot,E\setminus
A)>r_i\}$, with $r_i\downarrow 0$) provides the restriction to open
sets $[T]\res A$, satisfying $[T]\res A+[T]\res (E\setminus A)=[T]$
and ${\bf M}_p([T]\res A)\leq \|T\|_p(A)$. Since ${\bf M}_p$ is
subadditive and $[T]=[T]\res A+[T]\res C$, with $C=E\setminus A$, it
turns out that both inequalities are equalities:
\begin{equation}\label{localK}
{\bf M}_p([T]\res C)=\|T\|_p(C),\qquad {\bf M}_p([T]\res A)=
\|T\|_p(A).
\end{equation}
By \eqref{fasta} it follows also that
\begin{equation}\label{fastabis}
\fflat_p([T]\res\{u<s\}- [T]\res\{u<r\})\leq {\rm
Lip}(u)\|T\|_p(\{r\leq u<s\})
 \qquad\forall r,\,s\in\R,\,\,r<s,
\end{equation}
so that $r\mapsto [T]\res\{u<r\}$ is left continuous, as a map from
$\R$ to $\fcp{k}{E}$, and continuous out of a countable set
(contained in the set of $r$'s satisfying $\|T\|_p(\{u=r\})>0$).

\subsection{Cone construction}

Given $x\in E$ and $S\in\ic{k}{E}$, the cone construction in
\cite[Proposition~10.2]{ak2} provides a current in $\ic{k+1}{E}$,
that we shall denote by $\{x\}\times S$, supported on the union of
the segments joining $x$ to points in the support of $S$, and
satisfying
\begin{equation}\label{basicco}
\partial (\{x\}\times S)=S-\{x\}\times\partial S.
\end{equation}
In addition we have the inequality
\begin{equation}\label{basicco1}
\mass(\{x\}\times S)\leq r\mass(S)
\end{equation}
where $r$ is the radius of the smallest closed ball
$\overline{B}_r(x)$ containing the support of $S$. Since for
$R\in\ic{k}{E}$ and $S\in\ic{k+1}{E}$ we have
$$
\{x\}\times (R+\partial S)=\{x\}\times (R+S)-\partial (\{x\}\times
S),
$$
we immediately get $\fflat(\{x\}\times T)\leq 2{\rm
diam}(E)\fflat(T)$ for $T\in\ic{k}{E}$. By density and continuity
the cone construction uniquely extends to all $\fc{k}{E}$ and still
satisfies \eqref{basicco}. In addition, since $\ic{k}{E}$ is dense
in mass norm in $\rc{k}{E}$, and the approximation can easily be
done in such a way as to retain the bounds on the support (see
\cite{ambkatz}), we conclude that \eqref{basicco1} holds when
$S\in\rc{k}{E}$. In this case, it is proved in
\cite[Proposition~10.2]{ak2} that $\{x\}\times S\in\rc{k+1}{E}$, so
that
\begin{equation}\label{basicco3}
\fflat_p(\{x\}\times T)\leq 2{\rm diam}(E)\fflat_p(T).
\end{equation}

We will also need the inequality
\begin{equation}\label{basicco2}
\mass_p(\{x\}\times T)\leq r\mass_p(T)
\end{equation}
for all $T\in\fc{k}{E}$ with finite $\mass_p$ mass, whose measure
$\|T\|_p$ is supported in $\overline{B}_r(x)$. We sketch its simple
proof, based on \eqref{basicco1} and on the definition of $\mass_p$:
let $T_h\in\fc{k}{E}$ with $\mass(T_h)\to\mass_p(T)$ and
$\fflat_p(T_h-T)\to 0$ and $r'>r$. We know that
$\|T_h\|(\{d(\cdot,x)>(r+r')/2\})\to 0$, hence we can replace $T_h$
by its image $\tilde{T}_h$ under the 2-Lipschitz radial retraction
of $E$ onto the ball $\overline{B}_{(r+r')/2}(x)$ to obtain
$\tilde{T}_h$ supported on the ball, still converging to $T$ in
$\fflat_p$ distance and with $\mass(\tilde{T}_h)\to\mass_p(T)$. The
inequality \eqref{basicco3} yields the $\fflat_p$ convergence of
$\{x\}\times\tilde{T}_h$ to $\{x\}\times T$; then, passing to the
limit in \eqref{basicco1} gives $\mass_p(\{x\}\times T)\leq
\tfrac{1}{2}(r+r')\mass_p(T)$. Eventually we can let $r'\downarrow
r$ to obtain \eqref{basicco2}.

\subsection{Isoperimetric inequality}

The next result is proved in \cite[Corollary~8.7]{ambkatz}, adapting
the technique in \cite{We,We1}.

\begin{proposition}[Isoperimetric inequality in $\fcp{k}{E}$]
\label{pisop1} For $k\geq 1$ there exist constants $\delta_k$ such
that, if $[L]\in\fcp{k}{E}$ is a non zero current with $\partial
[L]=0$ and bounded support then
$$
\inf\left\{\frac{{\bf M}_p([T])}{\bigl[{\bf
M}_p([L])\bigr]^{(k+1)/k}}:\
[T]\in\fcp{k+1}{E},\,\,\partial[T]=[L]\right\}\leq\delta_k.
$$
\end{proposition}

\subsection{Slice operators}

Having defined the restriction to the sets $\{u<r\}$, $u\in {\rm
Lip}(E)$, the slice operator $T\in\fc{k}{E}\mapsto\langle
T,u,r\rangle\in\fc{k-1}{E}$ is defined by
$$
T\mapsto\langle T,u,r\rangle:= \partial (T\res\{u<r\})-(\partial
T)\res\{u<r\}
$$
whenever the right hand side makes sense (for $\Leb{1}$-a.e.
$r\in\R$). Since $\partial^2=0$ we have
$$
\partial \langle T,u,r\rangle=-\langle\partial T,u,r\rangle
\qquad\text{for $\Leb{1}$-a.e. $r\in\R$.}
$$
By \eqref{cachan} it follows that
\begin{equation}\label{cachaniv}
\int_m^{*\ell}\fflat(\langle T,u,r\rangle)\,dr\leq 2(\ell-m+{\rm
Lip}(u))\fflat(T)\qquad\forall m,\,\ell\in\R,\,\,m\leq\ell,
\end{equation}
\begin{equation}\label{cachanv}
\int_m^{*\ell}\fflat_p(\langle T,u,r\rangle)\,dr\leq 2(\ell-m+{\rm
Lip}(u))\fflat_p(T)\qquad\forall m,\,\ell\in\R,\,\,m\leq\ell.
\end{equation}

\subsection{0-dimensional chains}

It is not hard to show (see \cite[Proposition~8.4,
Theorem~8.5]{ambkatz} for a detailed proof) that, for
$T\in\fc{0}{E}$, we have the representation
\begin{equation}\label{dirac0}
\|T\|_p=\sum_{i=1}^N m_i\delta_{x_i}
\end{equation}
with $1\leq m_i\leq p/2$ and $x_i\in E$ distinct.

\subsection{Euclidean currents $\modp$ and flat chains with
coefficients in $\Z_p$}

In Euclidean spaces $\R^n$, a more general theory of currents with
coefficients in a normed Abelian group $G$ has been developed by
White in \cite{white1}, \cite{white2} on the basis of Fleming's
paper \cite{Fl}. The basic idea of \cite{Fl} is to consider the
abstract completion of the class of \emph{weakly polyhedral chains}
with respect to a flat distance. These objects are described by
finite sums
$$
\sum_i g_i \segop S_i\segcl,
$$
where $g_i\in G$ and $S_i$ are $k$-dimensional polyhedra, i.e. $S_i$
is contained in a $k$-plane and $\partial S_i$ is contained in
finitely many $(k-1)$-planes (we use the adjective \emph{weakly} to
avoid a potential confusion with the smaller class of polyhedral
currents of the deformation theorem: for these $S_i$ are $k$-cells
of a standard cubical decomposition of $\R^n$). The family of weakly
polyhedral chains with coefficients in $\Z_p$ has an obvious
additive structure inherited from $G$ and a boundary operator in
this class can be easily defined. The mass $\mass_G(T)$ of a weakly
polyhedral chain $T$ can be defined by minimizing $\sum_i
\|g_i\|\Haus{k}(S_i)$ among all possible decompositions of $T$, and
a flat distance is defined as follows:
$$
\fflat^P_G(T):=\inf\left\{\mass_G(R)+\mass_G(\partial S):\
\text{$R$, $S$ weakly polyhedral}\right\}.
$$
In the particular case $G=\Z_p$ we can obviously think weakly
polyhedral chains with coefficients in $\Z_p$ as currents with
coefficients in $\Z_p$ and the flat distance $\fflat_p^P=\fflat_G^P$
reads as follows:
\begin{equation}\label{flatp}
\fflat^P_p(T):=\inf\left\{\mass_p(R)+\mass_p(\partial S):\
\text{$R$, $S$ weakly polyhedral}\right\}.
\end{equation}
Obviously $\fflat^P_p(T)\geq\fflat_p(T)$ (because a larger class of
currents is considered in \eqref{may2}), but
Proposition~\ref{p1maggio} shows that the two flat distances are
equivalent in the class of weakly polyhedral chains. A direct
consequence of the equivalence of the two flat distances is the
following result, showing that currents $\modp$ are in canonical 1-1
correspondence with flat chains with coefficients in $\Z_p$ and that
the equivalence of the flat distances persists.

\begin{proposition}\label{pchess}
If $(T_i)$ is a Cauchy sequence with respect to $\fflat^P_p$, with
$T_i$ $k$-dimensional and weakly polyhedral, then
$\fflat_p(T_i-T)\to 0$ for some $T\in\fc{k}{E}$ with
$\fflat_p(T)\leq\lim_i\fflat_p^P(T_i)$. Conversely, if $T\in\fc{k}{E}$ then there exist $T_i$ $k$-dimensional and
weakly polyhedral with $T=\lim_i T_i$ with respect to $\fflat_p$; moreover $(T_i)$ is a Cauchy sequence with respect to
$\fflat_p^P$ and $\lim_i\fflat_p^P(T_i)\leq C\fflat_p(T)$. The
constant $C=C(n,k)$ is given by Proposition~\ref{p1maggio}.
\end{proposition}

\section{$\mass_p$ mass and slice mass}\label{s2}

In this section we introduce another notion of $p$ mass, the
so-called slice mass based on the 0-dimensional slices of the flat
chain, and we compare it with $\|T\|_p$.

We begin with some technical results stating more precise properties
of the slice operator. The first one concerns the inequality
\begin{equation}\label{mpslice}
\int_{\R}\|\langle T,\pi,r\rangle\|_p(E)\,dr\leq\|T\|_p(E)\qquad
\forall T\in\fc{k}{E}
\end{equation}
and all $\pi\in {\rm Lip}_1(E)$.
Let $(T_h)\subset\ic{k}{E}$ be satisfying
$\sum_h\fflat_p(T_h-T)<\infty$ and ${\bf M}(T_h)\to {\bf M}_p(T)$.
We know from \eqref{amb3} that $\sum_h\fflat_p(\partial T_h-\partial
T)<\infty$, hence
$$
\lim_{h\to\infty}T_h\res\{\pi<r\}=T\res\{\pi<r\},\qquad
\lim_{h\to\infty}(\partial T_h)\res\{\pi<r\}=(\partial
T)\res\{\pi<r\}
$$
with respect to the $\fflat_p$ distance for $\Leb{1}$-a.e. $r\in\R$.
It follows that $\fflat_p(\langle T_h-T,\pi,r\rangle)\to 0$ for
$\Leb{1}$-a.e. $r\in\R$.

First, let us check the measurability of $r\mapsto\|\langle
T,\pi,r\rangle\|_p(E)$. Since $\mass_p$ is lower semicontinuous in
$\fc{k}{E}$ we can find a nondecreasing sequence of
$\fflat_p$-continuous functions $G_i:\fc{k}{E}\to [0,+\infty)$ with
$G_i(T)\uparrow\mass_p(T)$ for all $T\in\fc{k}{E}$; taking into
account that
$$
\mass_p(\langle T,\pi,r\rangle)=\lim_{i\to\infty}\lim_{h\to\infty}
G_i(\langle T_h,\pi,r\rangle)
$$
for $\Leb{1}$-a.e. $r\in\R$, we need only to check the measurability
of $r\mapsto G_i(\langle S,\pi,r\rangle)$ for $S\in\rc{k}{E}$, which
is achieved in Lemma~\ref{ltechnical1}. The inequality
\eqref{mpslice} is known for $T\in\rc{k}{E}$ and for the ${\bf M}$
mass, see \cite[Theorem~5.7]{ak2}. Then, lower semicontinuity of
${\bf M}_p$ and Fatou's lemma give
\begin{eqnarray*}
\int_{\R}{\bf M}_p(\langle T,\pi,r\rangle)\,dr&\leq&
\int_{\R}\liminf_{h\to\infty}{\bf M}_p(\langle
T_h,\pi,r\rangle)\,dr\leq\liminf_{h\to\infty} \int_{\R}{\bf
M}_p(\langle T_h,\pi,r\rangle)\,dr\\&\leq&\liminf_{h\to\infty}
\int_{\R}{\bf M}(\langle
T_h,\pi,r\rangle)\,dr\leq\liminf_{h\to\infty} {\bf M}(T_h)={\bf
M}_p(T).
\end{eqnarray*}

\begin{lemma}[Slice and restriction commute]\label{lslice2}
Let $T\in\fc{k}{E}$, $\pi\in {\rm Lip}(E)$ and $u\in {\rm Lip}(E)$.
Then
\begin{equation}\label{immaco}
\langle T,\pi,r\rangle\res \{u<s\}=\langle T\res
\{u<s\},\pi,r\rangle \qquad\text{for $\Leb{2}$-a.e. $(r,s)\in\R^2$.}
\end{equation}
\end{lemma}

\noindent {\sc Proof.} The identity \eqref{immaco} is known when
$T\in\rc{k}{E}$. Indeed (see \cite[Theorem~5.7]{ak2}), the slices
$R_r$ of $R\in\rc{k}{E}$ are uniquely determined, up to Lebesgue
negligible sets, by the following two properties:

(a) $R_r$ is concentrated on $\pi^{-1}(r)$ for $\Leb{1}$-a.e.
$r\in\R$;

(b) $\int_{{\bf R}} \phi(r) R_r\,dr=R\res (\phi\circ\pi) d\pi$ for
all $\phi:\R\to\R$ bounded Borel.

It is then immediate to check that, for $s$ fixed, the currents in
the left hand side of \eqref{immaco} fulfil (a) and (b) with
$R=T\res \{u<s\}$, therefore they coincide with $\langle
R,\pi,r\rangle$ for $\Leb{1}$-a.e. $r\in\R$.

Let now $(T_h)\subset\rc{k}{E}$ with $\sum_h\fflat(T_h-T)<\infty$
and let us consider the identities
\begin{equation}\label{immaco1}
\langle T_h,\pi,r\rangle\res \{u<s\}=\langle T_h\res
\{u<s\},\pi,r\rangle \qquad\text{for $\Leb{2}$-a.e. $(r,s)\in\R^2$.}
\end{equation}
We know that $\sum_h\fflat(T_h\res\{u<s\}-T\res\{u<s\})<\infty$ for
$\Leb{1}$-a.e. $s\in\R$; for any $s$ for which this property holds,
we have that the right hand sides in \eqref{immaco1} converge to
$\langle T\res \{u<s\},\pi,r\rangle$ with respect to $\fflat$ for
$\Leb{1}$-a.e. $r\in\R$; on the other hand, we know also that
$\sum_h\fflat(\langle T_h,\pi,r\rangle-\langle
T,\pi,r\rangle)<\infty$ for $\Leb{1}$-a.e. $r\in\R$; for any $r$ for
which this property holds the left hand sides in \eqref{immaco1}
converge with respect to $\fflat$ to $\langle T_h,\pi,r\rangle\res
\{u<s\}$ for $\Leb{1}$-a.e. $s\in\R$.

Therefore, passing to the limit as $h\to\infty$ in \eqref{immaco1},
using Fubini's theorem, we conclude.\qed

We can now consider the local version of \eqref{mpslice}.

\begin{lemma}\label{lslice1}
For all $T\in\fc{k}{E}$, $\pi\in {\rm Lip}_1(E)$ and $B\subset E$
Borel the function $r\mapsto\|\langle T,\pi,r\rangle\|_p(B)$ is
Lebesgue measurable and
\begin{equation}\label{mpsliceter}
\int_{\R}\|\langle T,\pi,r\rangle)\|_p(B)\,dr\leq\|T\|_p(B).
\end{equation}
Furthermore, the support of $\|\langle T,\pi,r\rangle\|_p$ is
contained in $\pi^{-1}(r)\cap {\rm supp\,}\|T\|_p$ for
$\Leb{1}$-a.e. $r\in\R$.
\end{lemma}

\noindent {\sc Proof.} We consider a closed set $C\subset E$ and the
sets $C_s:=\{u<s\}$, $s>0$, where $u:= d(\cdot,C)$. Thanks to the commutativity of
slice and restriction, for $\Leb{1}$-a.e. $s>0$ we have $\langle
T,\pi,r\rangle\res \{u<s\}=\langle T\res \{u<s\},\pi,r\rangle$ for
$\Leb{1}$-a.e. $r\in\R$. Also, Fubini's theorem ensures that
$$
\|\langle T,\pi,r\rangle\|_p(C_s)=\mass_p(\langle T,\pi,r\rangle\res
\{u<s\})\qquad\text{for $\Leb{1}$-a.e. $r\in\R$}
$$
for $\Leb{1}$-a.e. $s>0$. Therefore, for any $s$ satisfying both
conditions we conclude that $\|\langle T,\pi,r\rangle\|_p(C_s)=
\mass_p(\langle T\res\{u<s\},\pi,r\rangle)$ for $\Leb{1}$-a.e.
$r\in\R$. Since we already proved that $r\mapsto\mass_p(\langle
T\res \{u<s\},\pi,r\rangle)$ is Lebesgue measurable, this proves
that the map $r\mapsto\|\langle T,\pi,r\rangle\|_p(C_s)$ is Lebesgue
measurable. Letting $s\downarrow 0$ the same is true for the map
$r\mapsto\|\langle T,\pi,r\rangle\|_p(C)$. The same argument allows
to prove \eqref{mpsliceter} from \eqref{mpslice}.

The class of Borel sets $B$ such that $r\mapsto \|\langle
T,\pi,r\rangle\|_p(B)$ is Lebesgue measurable contains the closed
sets and satisfies the stability assumptions of Dynkin's lemma,
therefore it coincides with the whole Borel $\sigma$-algebra.
Finally, by monotone approximation \eqref{mpsliceter} extends from
closed sets to open sets; if $B$ is Borel, by considering a
nonincreasing sequence of open sets $(A_h)$ such that
$\|T\|_p(A_h)\downarrow\|T\|_p(B)$ we extend the validity of
\eqref{mpsliceter} from open sets to Borel sets. Eventually,
choosing $A=E\setminus{\rm supp\,}\|T\|_p$ yields that $\|\langle
T,\pi,r\rangle\|_p(A)=0$ for $\Leb{1}$-a.e. $r\in\R$. \qed

Being defined on the whole of $\fc{k}{E}$ the slicing operator can
be obviously iterated, leading to the next definition.

\begin{definition}[Iterated slices]
For $T\in\fc{k}{E}$, $2\leq m\leq k$ and
$\pi=(\pi_1,\ldots,\pi_m)\in {\rm Lip}(E,\R^m)$ we define
recursively the slices $\langle T,\pi,x\rangle$, $x\in\R^m$, by
$$
\langle T,\pi,x\rangle:=\langle\langle
T,(\pi_1,\ldots,\pi_{m-1}),(x_1,\ldots,x_{m-1})\rangle,\pi_m,x_m\rangle.
$$
\end{definition}

Notice that the slices above are defined, as in the codimension 1
case, for $\Leb{m}$-a.e. $x\in\R^m$, they are given by
\begin{equation}\label{immaco4}
\langle T,\pi,x\rangle=\lim_{h\to\infty}\langle T_h,\pi,x\rangle
\qquad\text{in $\fc{k-m}{E}$}
\end{equation} whenever
$T_h\in\ic{k}{E}$ and $\sum_h\fflat(T_h-T)<\infty$ and the
definition is independent of $T_h$, up to $\Leb{k}$-negligible sets.
Moreover, a straightforward induction argument based on
Lemma~\ref{lslice2} gives
\begin{equation}\label{immacoiter}
\langle T,\pi,x\rangle\res \{u<s\}=\langle T\res
\{u<s\},\pi,x\rangle \qquad\text{for $\Leb{m+1}$-a.e.
$(x,s)\in\R^{m+1}$}
\end{equation}
for all $u\in {\rm Lip}(E)$.

Using \eqref{immaco4}, \eqref{immacoiter} and
Lemma~\ref{ltechnical1} as in the proof of \eqref{mpslice} and
Lemma~\ref{lslice1} we obtain:

\begin{lemma}\label{lslice3}
For all $T\in\fc{k}{E}$, $1\leq m\leq k$, $\pi\in \bigl[{\rm
Lip}_1(E)\bigr]^m$ and $B\subset E$ Borel the function
$x\mapsto\|\langle T,\pi,x\rangle\|_p(B)$ is Lebesgue measurable and
\begin{equation}\label{mpslicequater}
\int_{\R^m}\|\langle T,\pi,x\rangle\|_p(B)\,dx\leq\|T\|_p(B).
\end{equation}
Furthermore, the support of $\|\langle T,\pi,x\rangle\|_p$ is
contained in $\pi^{-1}(x)\cap {\rm supp\,}\|T\|_p$ for
$\Leb{m}$-a.e. $x\in\R^m$.
\end{lemma}

\begin{remark}[$BV$ regularity of slices]\label{rBV}{\rm
A direct consequence of \eqref{fasta} is that, for all
$T\in\fc{k}{E}$ with $\mass_p(T)$ finite, $s\mapsto T\res\{\pi<s\}$
has bounded variation in $\R\setminus N$ with respect to $\fflat_p$.
Since $N$ is Lebesgue negligible it follows that $T\res\{\pi<s\}$
has essential bounded variation, and its total variation measure
does not exceed ${\rm Lip}(\pi)\|T\|_p$. The same is true for the
slice map $r\mapsto \langle T,\pi,r\rangle$ of currents $T$ having
finite $\mass_p$ mass and boundary with finite $\mass_p$ mass, and
the total variation measure does not exceed ${\rm
Lip}(\pi)(\|T\|_p+\|\partial T\|_p)$.\\ For higher dimensional
slices, we can combine \eqref{mpslicequater} and the
characterization of metric $BV$ functions in terms of restrictions
to lines (see \cite[4.5.9]{federer} or \cite{afp} for the case of
real-valued maps and \cite{ambrosiomet} for the case of metric space
valued maps) to obtain that, for all $\pi\in \bigl[{\rm
Lip}(E)\bigr]^m$, $1\leq m\leq k$, we have $\langle
T,\pi,x\rangle\in MBV(\R^k,\fcp{k-m}{E})$ and its total variation
measure $\|D\langle T,\pi,x\rangle\|$ does not exceed
$$
\prod_{i=1}^m{\rm Lip}(\pi_i)\bigl(\|T\|_p+\|\partial T\|_p\bigr).
$$
}\end{remark}

Motivated by Lemma~\ref{lslice3}, for $1\leq m\leq k$, $\pi\in
\bigl[{\rm Lip}_1(E)\bigr]^m$ and $B\subset E$ Borel we define
\begin{equation}\label{realfl1}
\|T\res d\pi\|_p(B):=\int_{\R^m}\|\langle T,\pi,x\rangle\|_p(B)\,dx
\end{equation}
(the notation is reminiscent of the \emph{real} flat chain $T\res
d\pi=\int \langle T,\pi,x\rangle\,dx$). Notice that $\|T\res
d\pi\|_p$ is a $\sigma$-additive Borel measure less than $\|T\|_p$.

We shall also need the fact that $\|T\|_p$ has no atomic part:

\begin{lemma}\label{lnoatom}
The measure $\|T\|_p$ has no atom for all $T\in\fc{k}{E}$ with
finite $\mass_p$ mass.
\end{lemma}

\noindent {\sc Proof.} Writing $T=R+\partial S$ with $R\in\rc{k}{E}$
and $S\in\rc{k+1}{E}$, and noticing that
$\|T\|_p\leq\|R\|_p+\|\partial S\|_p\leq\|R\|+\|\partial S\|_p$,
since $\|R\|$ has no atom we can assume with no loss of generality
that $T=\partial S$. Fix $x\in E$, $\eps>0$ and $\bar r>0$ so small
that $\|S\|(B_{2\bar r}(x))<\eps$. Now, notice that
\begin{eqnarray*}
T\res\{d(x,\cdot)<s\}&=&\partial(S\res \{d(x,\cdot)<s\})- \langle
S,d(x,\cdot),s\rangle\\&=&
\partial (S\res\{d(x,\cdot)<s\}-\{x\}\times \langle
S,d(x,\cdot),s\rangle)+\{x\}\times\langle T,d(x,\cdot),s\rangle
\end{eqnarray*}
for $\Leb{1}$-a.e. $s>0$. Let now $r\leq \bar r$; since for $s<2r$
\eqref{basicco1} and \eqref{basicco2} give
$$
\mass(\{x\}\times \langle S,d(x,\cdot),s\rangle)\leq 2r\mass(\langle
S,d(x,\cdot),s\rangle),\quad \mass_p(\{x\}\times \langle
T,d(x,\cdot),s\rangle)\leq 2r\mass_p(\langle T,d(x,\cdot),s\rangle),
$$
we can choose $s\in (r,2r)$ and average to get
\begin{eqnarray*}
\frac{1}{r}\int_r^{*2r}\fflat_p(T\res\{d(x,\cdot)<s\})\,ds&\leq&\eps+
2\int_r^{2r}\mass(\langle S,d(x,\cdot),s\rangle)+\mass_p(\langle
T,d(x,\cdot),s\rangle)
\,ds\\&\leq&\eps+2\|S\|(B_{2r}(x)\setminus\{x\})+
2\|T\|_p(B_{2r}(x)\setminus\{x\})<2\eps
\end{eqnarray*}
for $r\leq \bar r$ small enough. Since $\eps>0$ is arbitrary, it
follows that we can find $(s_j)\downarrow 0$ such that
$T\res\{d(x,\cdot)<s_j\}\to 0$ $\modp$ and
$$
\mass_p(T-T\res\{d(x,\cdot)<s_j\})=\|T\|_p(\{d(x,\cdot)\geq
s_j\})\leq \mass_p(T)-\|T\|_p(\{x\}).
$$
Then, the lower semicontinuity of $\mass_p$ gives that
$\|T\|_p(\{x\})=0$. \qed

In the next theorem and in the sequel we will often deal with
exceptional sets depending on the slicing map $\pi$. For this reason
it will be convenient to restrict these maps to a sufficiently rich
but countable set: we fix a set ${\mathcal D}\subset {\rm Lip}_1(E)$
countable and dense in ${\rm Lip}_1(E)$ with respect to the sup
norm.

The next important result shows that currents with finite $\mass_p$
mass and boundary with finite $\mass_p$ mass are uniquely determined
by their $0$-dimensional slices. We don't know whether the result is
true for all flat chains with finite $\mass_p$ mass: we shall prove
this fact in a more restrictive class of spaces in Section~\ref{s6}.

\begin{theorem}\label{tabs1}
Let $T\in\fc{k}{E}$ with finite ${\bf M}_p$ mass and boundary
with finite ${\bf M}_p$ mass. Assume that, for
some $m\in [1,k]$ the following property holds:
$$ \text{for all $\pi\in [{\mathcal D}]^m$, $\langle
T,\pi,x\rangle=0$ $\modp$ for $\Leb{m}$-a.e. $x\in\R^m$.}$$ Then
$T=0$ $\modp$.
\end{theorem}
\noindent {\sc Proof.} We argue by induction on $m$ and we consider
first the case $m=1$. In the proof of the case $m=1$ we consider
first the case when $\partial T=0$ $\modp$, then the general case.

\noindent {\bf Step 1.} Assume $\partial T=0$ $\modp$. According to
the Lyapunov theorem the range of a finite nonnegative measure with
no atom is a closed interval. Hence, thanks to Lemma~\ref{lnoatom},
for any $\eps>0$ we can find a finite Borel partition
$B_1,\ldots,B_N$ of $E$ with $\|T\|_p(B_i)<\eps$; also, we can find
compact sets $K_i\subset B_i$ such that $\sum_i\|T\|_p(E\setminus
\cup_i K_i)\leq\eps$. Since the sets $K_i$ are pairwise disjoint, we
can find $\delta>0$ and $\phi_i\in {\mathcal D}$ such that, for
$r\in (\delta,2\delta)$, the open sets $A_i:=\{\phi_i<r\}$ are
pairwise disjoint, contain $K_i$ and satisfy $\|T\|_p(A_i)\leq\eps$
(just choose $\phi_i$ very close to $d(\cdot,K_i)$ and $2\delta$
less than the least distance between the $K_i$). In addition, for
$\Leb{1}$-a.e. $r\in (\delta,2\delta)$ the following property is
fulfilled:
$$
\partial (T\res A_i)=
\langle T,d(\cdot,K_i),r\rangle+(\partial T)\res A_i=0\,\,\,\modp
\qquad\text{for all $i=1,\ldots,N$.}
$$
Now we choose $r\in (\delta,2\delta)$ for which all the properties
above hold and we apply the isoperimetric inequality in $\fcp{k}{E}$
to obtain $S_i\in\fc{k+1}{E}$ with $\partial S_i=T\res A_i$ and
$${\bf M}_p(S_i)\leq\gamma_k\bigl({\bf M}_p(T\res
A_i)\bigr)^{1+1/k}\leq \eps^{1/k}\|T\|_p(A_i).$$ By applying the
cone construction to the cycle $T-\sum_iT\res A_i$ $\modp$, whose
${\bf M}_p$ mass is less than $\eps^{1/k}$, we obtain one more $S_0$
whose boundary $\modp$ is $T-\sum_iT\res A_i$ with mass less than
$2{\rm diam}(E)\eps$. It follows that
$$
\partial\sum_{i=0}^N [S_i]=[T], \qquad {\bf M}_p(\sum_{i=0}^NS_i)
\leq 2{\rm diam}(E)\eps+\gamma_k\eps^{1/k}{\bf M}_p(T).
$$
Since $\fflat_p(T)\leq\fflat_p(\sum_0^NS_i)\leq {\bf
M}_p(\sum_0^NS_i)$ and $\eps$ is arbitrary, this proves that
$[T]=0$.

\noindent {\bf Step 2.} The case $k=1$ is covered by
Corollary~\ref{cmain1d1} in Section~\ref{s3}: it shows the existence
of $T'\in\rc{1}{E}$ with $T'=T$ $\modp$, so that the slices of $T'$
vanish $\modp$ and therefore the multiplicity of $T'$ is $0$
$\modp$. In the case $k>1$ we can use the commutativity of slice and
restriction to show that the slices of $\partial T$ vanish $\modp$,
so that we can apply Step 1 to the cycle $\partial T$ to obtain
$\partial T=0$ $\modp$. Hence by applying Step 1 again we obtain
that $T=0$ $\modp$.

The proof of the induction step $m\mapsto m+1$ is not difficult: let
us fix $\pi\in\mathcal D$ and let us consider the codimension 1
slices $\langle T,\pi,t\rangle$, $\pi\in {\mathcal D}$; by
assumption, for all $q\in [{\mathcal D}]^m$, the $m$-codimensional
slices of $\langle T,\pi,t\rangle$ given by $\langle\langle
T,\pi,t\rangle,q,z\rangle=0$ vanish $\modp$ for $\Leb{m+1}$-a.e.
$(t,z)$; since ${\mathcal D}$ is countable we can find a
$\Leb{1}$-negligible set $N$ such that, for $t\notin N$ and for all
$q\in [{\mathcal D}]^m$ the slices vanish $\modp$ for $\Leb{m}$-a.e.
$z\in\R^m$. The induction assumption then gives $\langle
T,\pi,t\rangle=0$ $\modp$ for all $t\in\R\setminus N$. Eventually
the first step of the induction allows to conclude that $T=0$
$\modp$. \qed

\begin{definition}[Slice $\mass_p$ mass
$\|T\|_p^*$]\label{dslicemass} We define $\|T\|_p^*$ as the least
upper bound, in the lattice of nonnegative finite measures in $E$,
of the measures $\|T\res d\pi\|_p$, when $\pi$ runs in $\bigl[{\rm
Lip}_1(E)\bigr]^k$.
\end{definition}

Thanks to Theorem~\ref{tabs1} we know that $\|T\|_p^*$ provides a
reasonable notion of $p$-mass, since $\|T\|_p^*=0$ implies $T=0$
$\modp$, at least for flat chains $T$ whose finite $\mass_p$ mass
and boundary with finite $\mass_p$ mass. In addition,
\eqref{mpslicequater} with $m=k$ gives the inequality
$$
\|T\|_p^*\leq\|T\|_p,
$$
so that $\|T\|_p^*$ is well defined. We don't know, however, whether
equality holds in general, or whether (in case equality fails), an
isoperimetric inequality holds for $\|T\|_p^*$. In Section~\ref{s6} we
shall prove that the two notions of $p$-mass coincide in a suitable
class of spaces $E$.

\begin{corollary}\label{cabs3}
Let $T\in\fc{k}{E}$ with finite ${\bf M}_p$ mass. Then
$\|T\|_p^*(B)=0$ whenever $B$ is a $\Haus{k}$-negligible set.
\end{corollary}

\noindent {\sc Proof.} We fix $\pi\in \bigl[{\rm Lip}_1(E)\bigr]^k$
and we notice that, by the coarea inequality
\cite[Theorem~2.10.25]{federer}, $\Haus{0}(B\cap\pi^{-1}(x))=0$
(i.e. $B\cap\pi^{-1}(x)$ is empty) for $\Leb{k}$-a.e. $x\in\R^k$.
Also, Lemma~\ref{lslice3} gives that $\|\langle T,\pi,x\rangle\|_p$
is concentrated on $\pi^{-1}(x)$ for $\Leb{k}$-a.e. $x\in\R^k$.
Then, $\|\langle T,\pi,x\rangle\|_p(B)=0$ for $\Leb{k}$-a.e.
$x\in\R^k$, so that $\|T\res d\pi\|_p(B)=0$. Since $\pi$ is
arbitrary we conclude that $\|T\|_p^*(B)=0$.\qed

\section{Rectifiability in the case $k=1$}\label{s3}

Our goal in this section is to prove the rectifiability of
1-dimensional flat chains. We shall actually prove a more precise
version of Corollary~\ref{cmain3} when $\partial T=0$ $\modp$,
namely the existence of a cycle $T'\in\ic{1}{E}$ in the equivalence
class of $T$.

\begin{theorem}\label{tmain1d}
If $T\in\fc{1}{E}$ has finite $\mass_p$ mass and $\partial T=0$
$\modp$ then there exists $T'\in\ic{1}{E}$ with $\bdry T' = 0$ and
$T' = T$ $\modp$.
\end{theorem}

Writing any $T\in\fc{1}{E}$ with finite $\mass_p$ mass as
$R+\partial S$ with $R\in\rc{1}{E}$ and $S\in\rc{2}{E}$ we can apply
Theorem~\ref{tmain1d} to $\partial S$ to obtain the 1-dimensional
version of Corollary~\ref{cmain3}:

\begin{corollary}\label{cmain1d1}
For all $T\in\fc{1}{E}$ with finite $\mass_p$ mass there exists
$T'\in\rc{1}{E}$ with $T'=T$ $\modp$.
\end{corollary}

The proof of Theorem~\ref{tmain1d} follows by the construction of a
sequence $(T_n)\subset\ic{1}{E}$ of cycles satisfying
\begin{equation}\label{equation:good-sequence}
 \mass(T_n)\leq C\quad\text{and}\quad \fflat_p(T-T_n)\leq \frac{1}{n}
\end{equation}
for all $n\in\N$ and for a constant $C$ independent of $n$. Since
$E$ is a compact subset of a Banach space, we can then use the
closure and compactness theorems in \cite{ak2} to conclude that a
subsequence $(T_{n_j})$ converges weakly (i.e. in the duality with
all Lipschitz forms) to a cycle $T'\in\ic{1}{E}$. Since $E$ is
furthermore convex by \cite{Wenger-flatconv} we infer that $T_{n_j}$
converge in the flat norm to $T'$. It follows that $T =T'$ $\modp$
because
\begin{equation*}
 \fflat_p(T-T')\leq \fflat_p(T-T_{n_j}) + \fflat(T_{n_j}-T') \to 0
 \qquad\text{as $j\to\infty$.}
\end{equation*}

In order to construct a sequence $(T_n)$ of integral cycles satisfying
\eqref{equation:good-sequence} we proceed as follows. First we build, in
Lemma~\ref{weng1} below, approximating cycles $T_\varepsilon\in\ic{1}{E}$
whose boundary belongs to $p\cdot \ic{0}{E}$. Then, these cycles are in turn approximated
by finite sums $S$ of Lipschitz images of intervals, retaining the same boundary.
Eventually a combinatorial argument provides a cycle $S'$ in the same equivalence
class $\modp$ of $S$ with mass controlled by the mass of $S$.

\begin{lemma}\label{weng1}
 Let $[T]$ be as in the statement of Theorem~\ref{tmain1}. Then for every
$\varepsilon>0$ there exists $T_\varepsilon\in\ic{1}{E}$ such that
$\bdry T_\varepsilon\in p\cdot\ic{0}{E}$ and
\begin{equation*}
 \mass(T_\varepsilon)\leq \mass_p(T) + \varepsilon \quad\text{and}\quad
\fflat_p(T-T_\varepsilon)<\varepsilon.
 \end{equation*}
\end{lemma}

\noindent {\sc Proof.}
 Let $\varepsilon\in(0,1)$ and choose $T'\in\ic{1}{E}$ satisfying
$\fflat_p(T-T')<\varepsilon/3$ and
 \begin{equation*}
  \mass(T')\leq \mass_p(T) + \frac{\varepsilon}{3}.
 \end{equation*}
 Write $T=T' + R +\bdry S + p\cdot Q$ with $R\in\rc{1}{E}$, $S\in\rc{2}{E}$, $Q\in\ic{1}{E}$ and
 \begin{equation*}
  \mass(R) + \mass(S)\leq \frac{\varepsilon}{3}.
 \end{equation*}
 Since $\fflat_p(\bdry T) = 0$ we can write
$\bdry T = Z + \bdry U + p\cdot W$ with $Z\in\ic{0}{E}$, $U\in\rc{1}{E}$, $W\in\ic{0}{E}$ and
 \begin{equation*}
  \mass(Z) + \mass(U) \leq \frac{\varepsilon}{3}.
 \end{equation*}
 From this and the choice of $\varepsilon$ it follows that $Z=0$ and thus
 \begin{equation*}
  \bdry T' + \bdry R + p\cdot\bdry Q = \bdry U + p\cdot W.
 \end{equation*}
 Set $T_\varepsilon:= T' + R - U$. It is clear that $T_\varepsilon\in\rc{1}{E}$ and that
 \begin{equation*}
  \bdry T_\varepsilon = \bdry T' + \bdry R - \bdry U = p\cdot(W - \bdry Q) \in p\cdot\ic{0}{E},
 \end{equation*}
 so that $T_\varepsilon\in\ic{1}{E}$. Furthermore, we obtain
 \begin{equation*}
  \mass(T_\varepsilon) \leq \mass(T') + \mass(R) + \mass(U) \leq \mass_p(T) + \varepsilon
 \end{equation*}
 and
 \begin{equation*}
  T - T_\varepsilon = T - T' - R +U = \bdry S + U + p\cdot Q,
 \end{equation*}
 from which it follows that
 \begin{equation*}
  \fflat_p(T-T_\varepsilon) \leq \mass(U) + \mass(S) \leq \frac{2\varepsilon}{3}.
 \end{equation*}
 This concludes the proof.\qed

The following gives an almost optimal representation of currents in $\ic{1}{E}$ as
a superposition of curves. For a related result see
\cite{Wenger-sharpisop}, for the optimal result in $\R^n$ see \cite[4.2.25]{federer} (we shall
actually use this result in the proof).

\begin{lemma}\label{weng2}
Let $E$ be a length space and let $\tilde{T}\in\ic{1}{E}$. Then for every $\delta>0$ there exist
finitely many $(1+\delta)$-Lipschitz curves $c_i:[0,a_i]\to E$, $i=1,\ldots,n$, with image in
 $\overline{B}({\rm supp}\tilde{T}, \delta\mass(\tilde{T}))$
and such that $\bdry\tilde{T}=\sum (\Lbrack c_i(a_i)\Rbrack - \Lbrack c_i(0)\Rbrack)$,
 \begin{equation*}
  \mass(\bdry \tilde{T}) = \sum_{i=1}^n\mass(\Lbrack c_i(a_i)\Rbrack - \Lbrack c_i(0)\Rbrack),
 \end{equation*}
 \begin{equation*}
  \mass\left(\tilde{T}-\sum_{i=1}^n c_{i\sharp}\Lbrack\chi_{[0,a_i]}\Rbrack\right)\leq\delta\mass(\tilde{T}),
 \end{equation*}
 \begin{equation*}
  \sum_{i=1}^n a_i\leq (1+\delta)\mass(\tilde{T}).
\end{equation*}
\end{lemma}

\noindent {\sc Proof.}
Let $\delta'>0$ be small enough, to be determined later. Using Lemma 4 and Theorem 7 of \cite{kir1} one easily
shows that the existence of finitely many $(1+\delta')$-biLipschitz maps $\varphi_i: K_i\to E$, $i=1,\ldots, n$,
where the sets $K_i\subset\R$ are compact and such that $\varphi_i(K_i)\cap\varphi_j(K_j)=\emptyset$ if $i\not=j$, and
 \begin{equation}\label{equation:good-portion-T-covered}
  \|\tilde{T}\|\left(X\setminus\cup\varphi_i(K_i)\right)\leq\delta'\mass(\tilde{T}),
 \end{equation}
 see also \cite[Lemma 4.1]{ak2}.
 By McShane's extension theorem there exists a $(1+\delta')$-Lipschitz extension
$\overline{\eta}_i: E\to\R$ of $\varphi_i^{-1}$ for each $i=1,\ldots,n$.
 Now, write $\bdry\tilde{T}$ as $\bdry\tilde{T} = \sum_{i=1}^k (\Lbrack x_i\Rbrack -
 \Lbrack y_i\Rbrack)$ with $x_i\not=y_j$ for all $i, j$ (so that $\{x_1,\ldots,x_k\}$
 is the support of the positive part of $\partial\tilde{T}$ and
$\{y_1,\ldots,y_n\}$ is the support of the negative part).
 Note that $2k=\mass(\bdry\tilde{T})$.
 Set $\Omega:= \bigcup\varphi_i(K_i)\cup\{x_1, \ldots, x_k, y_1, \ldots, y_k\}$ and let
 $\{z_1,\ldots, z_m\}\subset\bigcup\varphi_i(K_i)$ be a finite and $\nu$-dense set for $\Omega$,
 where $\nu>0$ is such that
 \begin{equation}\label{equation:condition-K-distance}
  2\nu\leq \frac{\delta'}{1+\delta'}{\rm dist}(\varphi_i(K_i),\varphi_j(K_j))
 \quad\text{whenever $i\not=j$.}
 \end{equation}
 We set $N:= n+m+2k$ and define a map $\Psi: E\to \ell^\infty_N$ by
 \begin{equation*}
  \Psi(x):= \left(\overline{\eta}_1(x),\dots,\overline{\eta}_n(x), d(x,z_1),\dots,d(x,z_m), d(x, x_1),
d(x, y_1), \dots, d(x, x_k), d(x, y_k)\right).
 \end{equation*}
 Note that $\Psi$ is $(1+\delta')$-Lipschitz and $(1+\delta')$-biLipschitz on $\Omega$. Indeed, it is clear that
$\Psi$ is $(1+\delta')$-Lipschitz and that the restriction $\Psi\on{\varphi_i(K_i)}$ is $(1+\delta')$-biLipschitz
 for every $i$. Moreover, for $x\in \varphi_i(K_i)$ and $x'\in \varphi_j(K_j)$ with $i\not=j$ there exists
$z\in\varphi_i(K_i)$ with $d(x,z)\leq\nu$ and hence
 \begin{equation*}
  d(x,x')\leq d(x,z)+d(z,x')\leq 2 d(x,z)+d(z,x')-d(z,x)
 \leq\|\Psi(x')-\Psi(x)\|_\infty+2\nu,
 \end{equation*}
 from which the biLipschitz property on $\cup\varphi_i(K_i)$ follows together with
 \eqref{equation:condition-K-distance}. The other cases are trivial.
 By \cite[4.2.25]{federer} there exist at most countably many Lipschitz curves
 $\varrho_j:[0,a_j]\to\ell^\infty_N$ which are parametrized
 by arc-length, one-to-one on $[0,a_j)$ and which satisfy $\Psi_\sharp\tilde{T}= \sum_{j=1}^\infty
 \varrho_{j\sharp}\Lbrack\chi_{[0,a_j]}\Rbrack$ and
 \begin{equation}\label{equation:projected-length-curves}
  \mass(\Psi_\sharp\tilde{T})= \sum_{j=1}^\infty\mass(\varrho_{j\sharp}
  \Lbrack\chi_{[0,a_j]}\Rbrack)=\sum_{j=1}^\infty a_j
 \end{equation}
 and
 \begin{equation*}
 2k = \sum_{j=1}^\infty \mass(\Lbrack \varrho_j(a_j)\Rbrack - \Lbrack \varrho_j(0)\Rbrack).
 \end{equation*}
 After possibly reindexing the $\varrho_i$ and the $y_j$ we may assume without loss of generality that
 $\varrho_i(a_i) = \Psi(x_i)$ and $\varrho_i(0) = \Psi(y_i)$ for $i=1, \ldots, k$.
 It follows that $\varrho_i(a_i) = \varrho_i(0)$ for all $i\geq k+1$.
 Choose $M\geq k+1$ sufficiently large such that
 $R:= \sum_{j=M+1}^\infty \varrho_{j\sharp}\Lbrack\chi_{[0,a_j]}\Rbrack$ satisfies
 \begin{equation}\label{equation:R-small-endsum}
  \mass(R)\leq\delta'\mass(\tilde{T}).
 \end{equation}
 Since $E$ is a length space there exists a $(1+2\delta')$-Lipschitz extension $c_j:[0,a_j]\to E$ of
 $(\Psi\on{\Omega})^{-1}\circ(\varrho_j\on{\varrho_j^{-1}(\Psi(\Omega))})$ for each $j=1,\ldots, M$,
 and such that $c_j(a_j)=c_j(0)$ for $j=k+1, \ldots, M$.
 Note that $c_j(a_j) = x_j$ and $c_j(0)= y_j$ for all $j=1,\ldots, k$.
 We now have
  \begin{equation*}
  \sum_{j=1}^M \varrho_{j\#}\Lbrack\chi_{\varrho_j^{-1}(\Psi(\Omega)^c)}\Rbrack =
  \left[\Psi_\sharp(\tilde{T}\res\Omega^c)- R\right]\res\Psi(\Omega)^c,
 \end{equation*}
 from which it easily follows that
 \begin{equation*}
  \begin{split}
   T':&= \tilde{T}-\sum_{j=1}^M c_{j\sharp}\Lbrack\chi_{[0,a_j]}\Rbrack
   \\ &= \left(\Psi\on{\Omega}\right)^{-1}_\sharp\left[(R-\Psi_\sharp(\tilde{T}\res\Omega^c))\res\Psi(\Omega)\right]
   - \sum_{j=1}^M c_{j\sharp}\Lbrack\chi_{\varrho_j^{-1}(\Psi(\Omega)^c)}\Rbrack
   + \tilde{T}\res\Omega^c
  \end{split}
 \end{equation*}
 and, by moreover using \eqref{equation:good-portion-T-covered} and \eqref{equation:R-small-endsum},
 \begin{equation*}
  A:= \sum_{j=1}^M\Haus{1}(\varrho_j^{-1}(\Psi(\Omega)^c))
  =\sum_{j=1}^M\mass(\varrho_{j\sharp}\Lbrack\chi_{\varrho_j^{-1}(\Psi(\Omega)^c)}\Rbrack)
  \leq \delta'(2+\delta')\mass(\tilde{T}).
 \end{equation*}
 Using \eqref{equation:good-portion-T-covered}, \eqref{equation:R-small-endsum} and the facts that $c_i$ is
 $(1+2\delta')$-Lipschitz and $\Psi$ and $(\Psi\on{\Omega})^{-1}$ are $(1+\delta')$-Lipschitz, we obtain
 \begin{equation*}
 \begin{split}
  \mass(T')
   &\leq (1+\delta')\left[\mass(R) + (1+\delta')\|\tilde{T}\|(\Omega^c)\right] +(1+2\delta')A +
   \|\tilde{T}\|(\Omega^c)\\
   & \leq [5+8\delta'+3\delta'^2]\delta'\mass(\tilde{T}).
  \end{split}
 \end{equation*}
 Finally, using \eqref{equation:projected-length-curves} and the fact that $\Psi$ is $(1+\delta')$-Lipschitz,
 we estimate
 \begin{equation*}
  \sum_{j=1}^Ma_j\leq\mass(\Psi_\sharp\tilde{T})\leq (1+\delta')\mass(\tilde{T}).
 \end{equation*}
 This proves the statement given that $\delta'>0$ was chosen small enough.\qed

We now apply Lemma~\ref{weng2} to $\tilde{T}:=T_\varepsilon$, where $T_\varepsilon$ is
given by Lemma~\ref{weng1}. Set $T'':=\sum_{i=1}^n c_{i\sharp}\Lbrack\chi_{[0,a_i]}\Rbrack$.
We obtain, in particular,
\begin{equation}\label{equation:boundary-p-flat}
 \bdry T'' = \sum_{i=1}^n (\Lbrack c_i(a_i)\Rbrack - \Lbrack c_i(0)\Rbrack) \in p\cdot \ic{0}{E}.
\end{equation}

To conclude the proof of Theorem~\ref{tmain1d} we apply the following lemma.

\begin{lemma}\label{weng3}
 Let $E$ be a complete metric space, $n\geq 1$ and $p\geq 2$ integers. For each
 $i=1, \ldots, n$, let $c_i:[0,a_i]\to E$ be a Lipschitz curve. If
 $S:=\sum_{i=1}^n c_{i\sharp}\Lbrack\chi_{[0,a_i]}\Rbrack$ satisfies $\bdry S  \in p\cdot \ic{0}{E}$
then there exist $1\leq i_1<\cdots<i_k\leq n$ such that the current
\begin{equation}\label{equation:decomp-cycle-p}
 S':= S - p\cdot \sum_{j=1}^k c_{i_j\sharp}\Lbrack\chi_{[0,a_{i_j}]}\Rbrack
\end{equation}
is a cycle.
\end{lemma}

It follows in particular that $S - S' \in p\cdot\ic{1}{E}$ and that
\begin{equation*}
 \mass(S')\leq (p-1) \sum_{i=1}^n\mass(c_{i\sharp}\Lbrack\chi_{[0,a_i]}\Rbrack).
\end{equation*}

\noindent {\sc Proof.}
 It suffices to prove the lemma for the case that $c_i(a_i)\not=c_j(0)$ for all $1\leq i,j\leq n$, since we
can remove closed loops and we can concatenate $c_i$ and $c_j$
whenever $c_i(a_i)=c_j(0)$. Set $T_i:=
c_{i\sharp}\Lbrack\chi_{[0,a_i]}\Rbrack$.
 We first establish some notation: An ordered $k$-uple $(\alpha_0, \ldots,\alpha_k)$ with $k\geq 0$ and
 $\alpha_j\in\{1, \ldots, n\}$ is called admissible if either $k=0$ or $\alpha_r\not=\alpha_s$ for all
 $r\not=s$, $c_{\alpha_m}(a_{i_m}) = c_{\alpha_{m+1}}(a_{\alpha_{m+1}})$
 for all $m<k$ even, and $c_{\alpha_m}(0) = c_{\alpha_{m+1}}(0)$ for all $m<k$ odd.
 A decomposition $S = S_1 + \cdots + S_\ell + Q$ is called admissible if $Q\in p\cdot\ic{1}{E}$ and
 every $S_i$ is of the form
    \begin{equation*}
     S_i = \sum_{m=0}^{k_i} (-1)^m T_{\alpha(i, m)},
    \end{equation*}
 where $(\alpha(i,0),\ldots,\alpha(i,k_i))$ is an admissible $k_i$-uple for every $i\in [1,\ell]$,
 and if in addition the following properties hold:
   The sets  $$
 \Gamma_0:= \{\alpha(i,m):\ \text{ $m$ even, $1\leq i\leq\ell$}\},\qquad
 \Gamma_1:= \{\alpha(i,m):\ \text{ $m$ odd, $1\leq i\leq\ell$}\}$$
 are disjoint, $\Gamma_0\cup\Gamma_1= \{1, \ldots, n\}$, every index in $\Gamma_0$ appears exactly once,
 every index in $\Gamma_1$ appears exactly $p-1$ times, and $\bdry S_i=0$ if and only if $k_i$ is odd.
 Since the conditions imposed on admissible $k_i$-uples imply
 \begin{equation}\label{bdrysi}
 \partial S_i=
 \begin{cases}
  \segop c_{\alpha(i,k_i)}(a_{\alpha(i,k_i)})\segcl -\segop c_{\alpha(i,0)}(0)\segcl\neq 0 &
  \text{if $k_i$ is even}\\
  \segop c_{\alpha(i,k_i)}(0)\segcl -\segop c_{\alpha(i,0)}(0)\segcl &
  \text{if $k_i$ is odd}
 \end{cases}
 \end{equation}
 the last requirement is equivalent to the condition $c_{\alpha(i,k_i)}(0)=c_{\alpha(i,0)}(0)$
 whenever $k_i$ is odd.

  Note that, for example, the decomposition $S = S_1+\ldots+S_n$ with $S_i:= T_i$ is admissible
 ($Q=0$, $k_i=0$ for
  all $i$, $\Gamma_0=\{1,\ldots,n\}$, $\Gamma_1=\emptyset$). Suppose now that $S = S_1 + \cdots + S_\ell + Q$
  is an admissible decomposition. It is clear that if the set
 \begin{equation*}
  \Lambda(S_1,\ldots, S_\ell):= \{i\in\{1, \ldots, \ell\}: \bdry S_i\not=0\}
 \end{equation*}
 satisfies $|\Lambda|<p$ then in fact $\Lambda=\emptyset$. Indeed, $\sum_{i\in\Lambda}\bdry S_i=0$
 $\modp$ and \eqref{bdrysi} together with the fact that the indices $\alpha(i,k_i)$ appear only once
 (because $\bdry S_i\not=0$ implies $k_i$ even) give $\Lambda=\emptyset$.
 On the other hand, we claim that if $|\Lambda(S_1,\ldots, S_\ell)|\geq p$ then there exists an admissible
 decomposition $S = S_1' + \cdots + S_{\ell'}' + Q'$ with
\begin{equation*}
  |\Lambda(S'_1,\ldots, S'_{\ell'})| <  |\Lambda(S_1,\ldots, S_\ell)|.
\end{equation*}
In order to prove the claim, recall that $\bdry S_i= \Lbrack c_{\alpha(i,k_i)}(a_{\alpha(i,k_i)})\Rbrack -
\Lbrack c_{\alpha(i,0)}(0)\Rbrack$ whenever $k_i$ is even, and $\bdry S_i=0$ if $k_i$ is odd.
We call $\Lbrack c_{\alpha(i,k_i)}(a_{\alpha(i,k_i)})\Rbrack$ the right-boundary and
$\Lbrack c_{\alpha(i,0)}(0)\Rbrack$ the left-boundary of $S_i$.
We may now assume, up to a permutation of the $S_i$, that $\bdry
S_1\not=0$ and that the right-boundaries of  $S_2,\ldots,S_p$ equal
the right-boundary of $S_1$, with $\partial S_i\neq 0$ for
$i=2,\ldots,p$. We distinguish two cases:

First, suppose that $k_1=0$, so that $S_1=T_{\alpha(1,0)}$. Let
$r\in [1,p]$ be the number of integers $i\in [1,p]$ such that
$\partial S_i$ has the same left boundary of $\partial S_1$; we may
assume, again up to a permutation of the $S_2,\ldots,S_p$ and (in
case $r<p$) of $S_i$ for $i>p$ , that the left-boundaries of
$S_1,\ldots,S_r$ and (in the case $r<p$) $S_{p+1},\ldots,S_{2p-r}$
are all equal, with $\partial S_i\neq 0$ for $i=p+1,\ldots,2p-r$. We
define currents $S_j'$ by
\begin{equation*}
 S_j':= \left\{
 \begin{array}{l@{\quad}l}
   S_j - T_{\alpha(1,0)} & 2\leq j\leq r\\
   S_j - T_{\alpha(1,0)} + S_{p+j-r} & r+1\leq j \leq p.
  \end{array}\right.
\end{equation*}

Clearly, this yields an admissible decomposition $S =
S_2'+\cdots+S_p'+S_{2p-r+1}+\dots+S_\ell+Q'$
with $Q'=Q+pS_1$. Note that $S_j'$ is a cycle whenever $2\leq j\leq
r$ and therefore the number of non-cycles is strictly smaller if
$r\geq 2$; if $r=1$, since some non-cycles are concatenated in
groups of three, their total number is still strictly smaller in the
new decomposition.

Next, suppose that $k_1\geq 2$. Since $k_1$ is even the index
$\alpha(1,k_1)$ is in $\Gamma_0$ and thus appears exactly once.
Analogously, since $k_1-1$ is odd the index $\alpha(1,k_1-1)$ is in
$\Gamma_1$ and thus appears exactly $p-1$ times. We now construct a
new decomposition in which $\alpha(1,k_1)$ appears $p-1$ times and
$\alpha(1, k_1-1)$ only once. Let $r\in[1, p-1]$ be the number of
integers $i\in[1,p]$ such that $\alpha(i, t_i) = \alpha(1, k_1-1)$
for some $1\leq t_i\leq k_i$. Up to a permutation of the $S_2,
\dots, S_p$ and (in case $r<p-1$) of $S_i$ for $i>p$ we may
therefore assume that for every $i\in \{1,\ldots,
r\}\cup\{p+1,\ldots,2p-r-1\}$ we have $\alpha(i,t_i) =
\alpha(1,k_1-1)$ for a suitable $t_i$. If $i>p$ and $S_i$ is a cycle
then we may furthermore assume that $t_i=k_i$.
We now define the $S'_j$ as follows: First set
\begin{equation*}
 S'_1:= \sum_{m=0}^{k_1-2} (-1)^m T_{\alpha(1,m)}.
\end{equation*}
For $j\in\{2,\ldots,r\}$ we define a chain $S'_j$ and a cycle $S'_{l+j-1}$ such that
$S'_j+ S'_{l+j-1} = S_j - T_{\alpha(1, k_1)} + T_{\alpha(j, t_j)}$; more precisely,
let $S'_j$ be the `part' of $S_j$ preceding $\alpha(j, t_j)$ and $S'_{l+j-1}$ the
concatenation of the `part' of $S_j$ following $\alpha(j, t_j)$ with $-T_{\alpha(1, k_1)}$, thus
\begin{equation*}
 S_j':= \sum_{m=0}^{t_j-1} (-1)^m T_{\alpha(j,m)}\quad\text{ and }\quad
S_{\ell+j-1}':= \left[\sum_{m=t_j+1}^{k_j} (-1)^m T_{\alpha(j, m)}\right] - T_{\alpha(1, k_1)}.
\end{equation*}
Since the left- and right-boundaries of $T_{\alpha(1, k_1)}$ and the term in brackets in the above
equation agree, $S'_{\ell+j-1}$ is a cycle.
Let now $j\in\{r+1,\ldots,p-1\}$ and observe that the right-boundaries of $S_1$ and $S_j$ agree. Define
\begin{equation*}
S_j':=  S_j - T_{\alpha(1,k_1)}+\sum_{m=\bar{t}_j+1}^{\bar{k}_j} (-1)^m T_{\alpha(p-r+j, m)}\quad\text{ and }\quad
 S'_{p-r+j}:= \sum_{m=0}^{\bar{t}_j-1} (-1)^m T_{\alpha(p-r+j,m)},
\end{equation*}
where $\bar{t}_j = -1$ and $\bar{k}_j=k_{p-r+j}-1$ if $\bdry S_{p-r+j}=0$ and $\bar{t}_j=t_{p-r+j}$ and
$\bar{k}_j=k_{p-r+j}$ otherwise.
In particular, if $\bdry S_{p-r+j}=0$ then $S'_{p-r+j}=0$. Finally, set
$S_p':= S_p - T_{\alpha(1,k_1)} + T_{\alpha(1,k_1-1)}$, and for $j\in\{2p-r, \dots, \ell\}$ set $S'_j:= S_j$.
Observe that the index $\alpha(1, k_1)$ appears exactly $p-1$ times, $\alpha(1, k_1-1)$ exactly once,
and that all other indices appear exactly the same number of times as in the original admissible decomposition.
We therefore obtain an admissible decomposition $S = S'_1+\dots+S'_{\ell+r-1}+Q$.
Note that it has same number of chains with
non-zero boundary, however $S_1'$ has two edges less than $S_1$, that is, $k_1'=k_1-2$.
Applying the same procedure finitely many times allows us to reduce the second case to the first one.
This completes the proof of the claim.

We can now apply this claim finitely many times to obtain an admissible decomposition
$S = S_1 + \cdots + S_\ell + Q$ in which all $S_i$ are cycles. The current $S':= S_1+\cdots+S_\ell$
is clearly of the form (\ref{equation:decomp-cycle-p}) because
$$
S-\sum_{i=1}^\ell S_i=\sum_{i\in\Gamma_1}T_i-
\sum_{i=1}^\ell (-1)^m\sum_{\text{$m=1$, $m$ odd}}^{k_i}T_{\alpha(i,m)}=
p\sum_{i\in\Gamma_1}T_i
$$
and thus the proof of the lemma is complete with $\{i_1,\ldots,i_k\}=\Gamma_1$.\qed

\section{Lusin approximation of Borel maps by Lipschitz
maps}\label{s4}

Let $f:A\subset\R^k\to E$ a Borel map. For $x\in A$ we define
$\delta_xf$ as the smallest $M\geq 0$ such that
$$
\lim_{r\downarrow 0}r^{-k}\Leb{k}\bigl(\{y\in A\cap B_r(x):\
\frac{d(f(y),f(x))}{r}>M\}\bigr)=0.
$$
This definition is a simplified version of Federer's definition of
approximate upper limit of the difference quotients (we replaced
$|y-x|$ by $r$ in the denominator), but sufficient for our purposes.

\begin{theorem}\label{tdifftotale}
Let $f:A\subset\R^k\to E$ be Borel.
\begin{itemize}
\item[(i)] Let $k=n+m$, $x=(z,y)$, and assume that there exist Borel subsets
$A_1,\,A_2$ of $A$ such that $\delta_z(f(\cdot,y))<\infty$ for all
$(z,y)\in A_1$ and $\delta_y(f(z,\cdot))<\infty$ for all $(z,y)\in
A_2$. Then $\delta_x f<\infty$ for $\Leb{k}$-a.e. $x\in A_1\cap
A_2$;
\item[(ii)] if $\delta_xf<\infty$ for $\Leb{k}$-a.e. $x\in A$
there exists a sequence of Borel sets $B_n\subset A$ such that
$\Leb{k}(A\setminus\cup_n B_n)=0$ and the restriction of $f$ to
$B_n$ is Lipschitz for all $n$.
\end{itemize}
\end{theorem}

\noindent {\sc Proof.} For real-valued maps this result is basically
contained in \cite[Theorem~3.1.4]{federer}, with slightly different
definitions: here, to simplify matters as much as possible,
we avoid to mention any differentiability result.

\noindent (i) By an exhaustion argument we can assume with no loss
of generality that, for some constant $N$, $\delta_z f< N$ in $A_1$
and $\delta_y f< N$ in $A_2$. Moreover, by Egorov theorem (which
allows to transform pointwise limits, in our case as $r\downarrow
0$, into uniform ones, at the expense of passing to a slightly
smaller domain in measure), we can also assume that
\begin{equation}\label{egorov1}
\lim_{r\downarrow 0} r^{-m}\Leb{m}\bigl(\{y'\in B^m_r(y):\
\frac{d(f(z,y'),f(z,y))}{r}>N\}\bigr)=0 \qquad\text{uniformly for
$(z,y)\in A_2$.}
\end{equation}
We are going to show that $\delta_x f\leq 2N$ $\Leb{k}$-a.e. in
$A_1\cap A_2$. By the triangle inequality, it suffices to show that
\begin{equation}\label{rade1}
\lim_{r\downarrow 0} r^{-k}\Leb{k}\bigl(\{(z',y')\in B_r((z,y)):\
\frac{d(f(z',y),f(z,y))}{r}>N\}\bigr)=0
\end{equation}
and
\begin{equation}\label{rade2}
\lim_{r\downarrow 0} r^{-k}\Leb{k}\bigl(\{(z',y')\in B_r((z,y)):\
\frac{d(f(z',y'),f(z',y))}{r}>N\}\bigr)=0
\end{equation}
for $\Leb{k}$-a.e. $(z,y)\in A_1\cap A_2$. The first property is
clearly satisfied at all $(z,y)\in A_1$, because the sets in
\eqref{rade1} are contained in
$$
\{z'\in B_r(z):\ \frac{d(f(z',y),f(z,y))}{r}>N\}\times B_r(y).
$$
In order to show the second property \eqref{rade2} we can estimate
the quantity therein by
$$
r^{-n}\int_{B^1_r(z)}r^{-m}\Leb{m}\bigl(\{y'\in B_r(y):\
\frac{d(f(z',y'),f(z',y)}{r}>M\}\bigr)\,dz'+
\frac{\Leb{m}(B_r(y))\Leb{n}(B_r^2(z))}{r^{m+n}},
$$
where $B^1_r(z)=\{z'\in B_r(z):\ (z',y)\in A_2\}$ and
$B^2_r(z):=B_r(z)\setminus B^1_r(z)$. If we let $r\downarrow 0$, the
first term gives no contribution thanks to \eqref{egorov1}; the
second one gives no contribution as well provided that $z$ is a
density point in $\R^n$ for the slice $(A_2)_y:=\{z':\ (y,z')\in
A_2\}$. Since, for all $y$, $\Leb{n}$-a.e. point of $(A_2)_y$ is a
density point $(A_2)_y$, by Fubini's theorem we get that
$\Leb{k}$-a.e. $(y,z)\in A_2$ has this property.

\noindent (ii) Let $e_0\in E$ be fixed. Denote by $C_N$ the subset
of $A$ where both $\delta_x$ and $d(f,e_0)$ do not exceed $N$. Since
the union of $C_N$ covers $\Leb{k}$-almost all of $A$, it suffices
to find a family $(B_n)$ with the required properties covering
$\Leb{k}$-almost all of $C_N$. Let $\chi_k$ be a geometric constant
defined by the property
$$
\Leb{k}(B_{|x_1-x_2|}(x_1)\cap B_{|x_1-x_2|}(x_2))=\chi_k
\Leb{k}(B_{|x_1-x_2|}(0)).
$$
We choose $B_n\subset C_N$ and $r_n>0$ in such a way that
$\Leb{k}(C_N\setminus\cup_n B_n)=0$ and, for all $x\in B_n$ and
$r\in (0,r_n)$, we have
\begin{equation}\label{egorov}
\Leb{k}\bigl(\{y\in B_r(x):\
\frac{d(f(y),f(x))}{r}>N+1\}\bigr)\leq\frac{\chi_k}{2}\Leb{k}(B_r(x)).
\end{equation}
The existence of $B_n$ is again ensured by Egorov theorem.

We now claim that the restriction of $f$ to $C_n$ is Lipschitz.
Indeed, take $x_1,\,x_2\in B_n$: if $|x_1-x_2|\geq r_n$ we estimate
$d(f(x_1),f(x_2))$ simply with $4 r_n^{-1}\sup_{B_n}d(f,e_0)
|x_1-x_2|$. If not, by \eqref{egorov} at $x=x_i$ with $r=|x_1-x_2|$
and our choice of $\chi_k$ we can find $y\in B_r(x_1)\cap B_r(x_2)$
where
$$
\frac{d(f(y),f(x_1))}{r}\leq N+1\quad\text{and}\quad
\frac{d(f(y),f(x_2))}{r}\leq N+1.
$$
It follows that $d(f(x_1),f(x_2))\leq 2(N+1)|x_1-x_2|$.\qed

\begin{proposition}\label{pdeltametric}
Let $K\subset \Gamma\subset E$, with $K$ countably
$\Haus{1}$-rectifiable, and let $\pi\in {\rm Lip}(E)$ be injective
on $\Gamma$. Then $\delta (\pi\vert_\Gamma)^{-1}$ is finite
$\Leb{1}$-a.e. on $\pi(K)$.
\end{proposition}

\noindent {\sc Proof.} Assume first that $K=f(C)$ with $C\subset\R$
closed and $f:C\to K$ Lipschitz and invertible. The condition
$\delta (\pi\vert_\Gamma)^{-1}<\infty$ clearly holds at all points
$t=\pi(x)$ of density 1 for $\pi(K)$, with $x\in K$ satisfying
$$
\liminf_{y\in K\to x}\frac{|\pi(y)-\pi(x)|}{d(y,x)}>0.
$$
Indeed, at these points $t=\pi(x)$ we have
$x=(\pi\vert_\Gamma)^{-1}(t)$ and
$$
\liminf_{s\in \pi(K)\to
t}\frac{|(\pi\vert_\Gamma)^{-1}(s)-x|}{|s-t|}<\infty.
$$
If $N\subset K$ is the set where the condition above fails, the
Lipschitz function $p=\pi\circ f$ has null derivative at all points
in $f^{-1}(N)$ where it is differentiable, hence
$\Leb{1}(p(f^{-1}(N)))=0$. It follows that $\Leb{1}(\pi(N))=0$.

In the general case, write $K=N\cup\cup_i K_i$ with $\Haus{1}(N)=0$
and $K_i=f_i(C_i)$ pairwise disjoint, with $C_i\subset\R$ closed and
$f_i:C_i\to K_i$ Lipschitz and invertible. Let $B_i\subset\pi(K_i)$
be Borel sets such that the inverse $g_i$ of $\pi\vert_{K_i}$
satisfies $\delta g_i<\infty$ on $B_i$ and
$\Leb{1}(\pi(K_i)\setminus B_i)=0$. Since $\Haus{1}(\pi(N))=0$, the
union $\cup_i\pi(K_i)$ covers $\Leb{1}$-almost all of $\pi(K)$.
Hence, it suffices to show that $\delta (\pi\vert_\Gamma)^{-1}<\infty$ at
all points of density 1 for one of the sets $B_i$. This property
easily follows from the definition of $\delta$ and from the fact
that $(\pi\vert_\Gamma)^{-1}$ and $g_i$ coincide on $B_i$.\qed

\section{Countable rectifiability of $\|T\|_p^*$ in the case $k>1$}\label{s5}

In this section we show that the slice mass $\|T\|_p^*$ is
concentrated on a countably $\Haus{k}$-rectifiable set, adapting to
this context White's argument \cite{white1}; this provides a first
step towards the proof of Theorem~\ref{tmain1}.

The next technical lemma provides a useful commutativity property of
the iterated slice operator.

\begin{lemma}[Commutativity of slices]\label{lcommu}
Let $T\in\fc{k}{E}$ and $\pi=(p,q)$ with $p\in {\rm
Lip}(E;\R^{m_1})$, $q\in {\rm Lip}(E;\R^{m_2})$, $m_i\geq 1$ and
$m_1+m_2\leq k$. Then
\begin{equation}\label{alforno}
\langle \langle T,p,z\rangle,q,y\rangle=(-1)^{m_1m_2}\langle \langle
T,q,y\rangle,p,z\rangle \qquad\text{for $\Leb{m_1+m_2}$-a.e.
$(z,y)\in\R^{m_1}\times\R^{m_2}$.}
\end{equation}
\end{lemma}

\noindent {\sc Proof.} If $T\in\rc{k}{E}$ we know by
\cite[Theorem~5.7]{ak2} that the slices $S_{yz}=\langle\langle
T,p,z\rangle,q,y\rangle$ are characterized by the following two
properties:
\begin{itemize}
\item[(a)] $S_{yz}$ is concentrated on $p^{-1}(z)\cap q^{-1}(y)$ for
$\Leb{m_1+m_2}$-a.e. $(z,y)$;
\item[(b)] $\int \psi(y,z) S_{yz}\,dydz=T\res \psi(p,q)dp\wedge dq$
as $(k-m_1-m_2)$-dimensional currents for all bounded Borel
functions $\psi$.
\end{itemize}
It is immediate to check that $\langle \langle
T,q,y\rangle,p,z\rangle$ satisfy (a) and
$$
\int \psi(y,z)\langle \langle T,q,y\rangle,p,z\rangle\,dydz= T\res
\psi(p,q)dq\wedge dp=(-1)^{m_1+m_2} T\res\psi(p,q) dp\wedge dq,
$$
hence \eqref{alforno} holds. The general case can be achieved using
\eqref{immaco4}, choosing a sequence $(T_h)\subset\ic{k}{E}$ with
$\sum_h\fflat(T_h-T)<\infty$.\qed

In the next proposition we consider first the rectifiability of the
measures $\|T\res d\pi\|_p$ for $\pi$ fixed.

\begin{proposition}\label{prettipi}
Let $T\in\fc{k}{E}$ with finite ${\bf M}_p$ mass. Then, for all
$\pi\in \bigl[{\rm Lip}_1(E)\bigr]^k$, $\|T\res d\pi\|_p$ is
concentrated on a countably $\Haus{k}$-rectifiable set.
\end{proposition}

\noindent {\sc Proof.} By \eqref{dirac0} we obtain that $\|\langle
T,\pi,x\rangle\|_p$ consists for $\Leb{k}$-a.e. $x$ of a finite sum
of Dirac masses, with weights between 1 and $p/2$. Hence, we can
define
$$
\Lambda(x):=\left\{y\in\pi^{-1}(x):\ y\in {\rm supp\,}\|\langle
T,\pi,x\rangle\|_p\right\}
$$
and we can check that the set-valued function $\Lambda$ fulfils the
measurability assumption of Lemma~\ref{ltechnical2}. Indeed, for all
Borel sets $B$
$$
\left\{x\in\R^k:\ \Lambda(x)\cap B\neq\emptyset\right\}=
\left\{x\in\R^k:\ \|\langle T,\pi,x\rangle\|_p(B)>0\right\}
$$
and we know that the latter set is measurable, thanks to
Lemma~\ref{lslice3}. By Lemma~\ref{ltechnical2} we obtain disjoint
measurable sets $B_n=\{x: {\rm card\,} \Lambda(x)=n\}$ and measurable
maps $f_{j_1},\ldots,f_{j_n}$ satisfying \eqref{reponbn}.

Obviously it suffices to show that, for $n$ fixed and $C\subset B_n$
compact, the measure $B\mapsto\int_C\|\langle
T,\pi,x\rangle\|_p(B)\,dx$ is concentrated on a countably
$\Haus{k}$-rectifiable set. By a further approximation based on
Lusin's theorem we can also assume that $f_{j_1},\ldots,f_{j_n}$ are
continuous in $C$. Finally, since $f_{j_i}(x)\neq f_{j_\ell}(x)$
whenever $x\in B_n$ and $i\neq\ell$ we can also assume that the sets
$K_i:=f_{j_i}(C)$, $i=1,\ldots,n$, are pairwise disjoint. Notice that
$\pi:K_i\to C$ is injective and its inverse is $f_{j_i}$.

We consider now $u_i=d(\cdot,K_i)$ and let $s>0$ be the least
distance between the sets $K_i$, so that for $s_i\in (0,s/2)$ the
sets $\{u_i<s\}$ are pairwise disjoint; thanks to the commutativity
of slice and restriction, for $\Leb{1}$-a.e $s_i>0$ we have
\begin{equation}\label{armadio}
\langle T\res\{u_i<s_i\},\pi,x\rangle=\langle
T,\pi,x\rangle\res\{u_i<s_i\} \qquad\text{for $\Leb{k}$-a.e. $x$}
\end{equation}
for $i=1,\ldots,n$. Choosing $s_i\in (0,s/2)$ with this property and
setting $T_i:=T\res\{u_i<s_i\}$, we have
$$
\int_C\|\langle T,\pi,x\rangle\|_p(B)\,dx=\sum_{i=1}^n
\int_C\|\langle T_i,\pi,x\rangle\|_p(B)\,dx
$$
and it suffices to show that all measures $\mu_i(B):=\int_C\|\langle
T_i,\pi,x\rangle\|_p(B)\,dx$ are concentrated on a countably
$\Haus{k}$-rectifiable set. By \eqref{armadio} it follows that
$\Leb{k}$-almost all measures $\|\langle T_i,\pi,x\rangle\|_p$,
$x\in C$, are Dirac masses concentrated on $f_{j_i}(x)$.

We now fix $i$ and prove that $\mu_i$ is concentrated on a countably
$\Haus{k}$-rectifiable set by applying Theorem~\ref{tdifftotale}(ii)
to the inverse $f_{j_i}$ of $\pi\vert_{K_i}$. Let us consider the
sets
$$
C_z:=\left\{t\in\R:\ (z,t)\in C\right\},\qquad K_{iz}:=\left\{x\in
K_i:\ (\pi_1,\ldots,\pi_{n-1})(x)=z\right\}
$$
and the maps $g_z(t):=f_{j_i}(z,t):C_z\to K_{iz}$. We claim that,
for $\Leb{k-1}$-a.e. $z$, $\delta_t g_z<\infty$ $\Leb{1}$-a.e. in
$C_z$. Indeed, writing $x=(z,t)$ and
$$
\langle T_i,\pi,x\rangle=\langle S_z,\pi_k,t\rangle\qquad
\text{with}\qquad S_z:=\langle T_i,(\pi_1,\ldots,\pi_{k-1}),z\rangle
$$
we know that for $\Leb{k-1}$-a.e. $z$ the flat chain
$S_z\in\fc{1}{E}$ has finite $\mass_p$ mass and $\|\langle
S_z,\pi_k,t\rangle\|_p$ is a Dirac mass on $g_{iz}(t)$ for
$\Leb{1}$-a.e. $t\in C_z$.

We fix now $z$ with these properties; since $S_z\in\fc{1}{E}$,
Corollary~\ref{cmain1d1} provides $S_z'\in\rc{1}{E}$ with $S_z'=S_z$
$\modp$. Then, we can find countably $\Haus{1}$-rectifiable set $G$
on which $S_z'$ is concentrated, and therefore a
$\Leb{1}$-negligible set $N_z$ such that $\langle
S_z',\pi_k,t\rangle$ is concentrated on $G$ for all $t\in\R\setminus
N_z$. But since $\langle S_z',\pi_k,t\rangle=\langle
S_z,\pi_k,t\rangle$ $\modp$ for $\Leb{1}$-a.e. $t$, possibly adding
to $N_z$ another $\Leb{1}$-negligible set we can assume that
$\|\langle S_z',\pi_k,t\rangle\|_p$ is a Dirac mass on $g_{iz}(t)\in
G$ for all $t\in C_z\setminus N_z$. We denote by
$\tilde{K}_{iz}\subset G$ the set
$$
\tilde{K}_{iz}:=\left\{g_{iz}(t):\ t\in C_z\setminus N_z\right\}
$$
which is countably $\Haus{1}$-rectifiable as well and contained in
$K_{iz}$. Notice also that
$\Leb{1}(\pi_k(K_{iz}\setminus\tilde{K}_{iz}))=0$ because this set
is contained in $N_z$. Since $\pi_k\vert_{K_{iz}}$ is injective, we
can now apply Proposition~\ref{pdeltametric} with $K=\tilde{K}_{iz}$
and $\Gamma=K_{iz}$ to obtain that
$\delta((\pi_k)\vert_{K_{iz}})^{-1}<\infty$ $\Leb{1}$-a.e. on
$\pi_k(\tilde{K}_{iz})$ and therefore $\Leb{1}$-a.e. on
$\pi_k(K_{iz})$. But, since the inverse of $\pi\vert_{K_i}$ is
$f_{j_i}$, the inverse of $(\pi_k)\vert_{K_{iz}}$ is $g_z$. It
follows that $\delta g_z<\infty$ $\Leb{1}$-a.e. on $C_z$.

This proves the claim. Thanks to the commutativity of slice and
restriction, a similar property is fulfilled by $f_{j_i}$ with
respect to the other $(k-1)$ variables, hence
Theorem~\ref{tdifftotale}(i) ensures that $\delta_x f_{j_i}<\infty$
$\Leb{k}$-a.e. on $C$. This ensures that
Theorem~\ref{tdifftotale}(ii) is applicable to $f_{j_i}$, and in
turn the fact that $\mu_i$ is concentrated on a
$\Haus{k}$-rectifiable set.\qed

We recall that the supremum ${\mathcal M}-\sup_{i\in I}\mu_i$ of a
family of measures $\{\mu_i\}_{i\in I}$ is the smallest measure
greater than all $\mu_i$; it can be constructively defined by
\begin{equation}\label{explicit}
{\mathcal M}-\sup_{i\in I}\mu_i (B):= \sup \sum_{j=1}^N\mu_{i_j}(B_j)
\end{equation}
where the supremum runs among all finite Borel partitions
$B_1,\ldots,B_N$ of $B$, with $i_1,\ldots,i_N\in I$.

\begin{proposition}\label{psabato}
Let $T\in\fc{k}{E}$ with finite $\mass_p$ mass. Then $\|T\|_p^*$ is
concentrated on a countably $\Haus{k}$-rectifiable set.
\end{proposition}

\noindent {\sc Proof.} Let $I$ be an index set for $\bigl[{\rm
Lip}_1(E)]^k$, and consider for any $n\in\N$ a finite set
$J_n\subset I$ such that
$$\|T\|_p^*(E)\leq{\mathcal M}-\sup_{i\in J_n}\|T\res
d\pi_i\|_p(E)+2^{-n}$$
(its existence is a direct consequence of
\eqref{explicit}). Then, denoting by $J$ the union of the sets
$J_n$, the measure
$$
\sigma:={\mathcal M}-\sup_{i\in J}\|T\res d\pi_i\|_p
$$
is smaller than $\|T\|_p^*$ and with the same total mass, hence it
coincides with $\|T\|_p^*$. Since $J$ is countable, a countably
$\Haus{k}$-rectifiable concentration set for $\sigma$ can be
obtained by taking the union of countably $\Haus{k}$-rectifiable
sets, given by Proposition~\ref{prettipi}, on which the measures
$\|T\res d\pi_i\|_p$, $i\in J$, are concentrated.\qed

\section{Absolute continuity of $\|T\|_p$}\label{s6}

In this section we prove the absolute continuity of $\|T\|_p$ with
respect to $\|T\|_p^*$, and therefore the fact that also $\|T\|_p$
is concentrated on a countably $\Haus{k}$-rectifiable set. Then, we
can prove, using the isoperimetric inequality, density lower bounds
for $\|T\|_p$; these imply that the (minimal) concentration set has
actually finite $\Haus{k}$-measure.

The absolute continuity of $\|T\|_p$ depends on the following
extension of Theorem~\ref{tabs1} to all flat chains with finite
$\mass_p$ mass. We are presently able to prove this extension,
relying on the finite-dimensional results in \cite{white1} (in turn
based on the deformation theorem in \cite{white2}), only in a
smaller class of spaces $E$.

\begin{definition}\label{dsfda}
We say that a Banach space $(F,\|\cdot\|)$ has the strong
finite-dimensional approximation property if there exist maps
$\pi_n:F\to F$, with uniformly bounded Lipschitz constants, such
that $\pi_n(F)$ is contained in a finite-dimensional subspace of $F$
and $$\lim_{n\to\infty}\|x-\pi_n(x)\|=0\qquad\forall x\in F.$$
\end{definition}

Obviously all spaces having a Schauder basis (and, in particular,
separable Hilbert spaces) have the strong finite-dimensional
approximation property and, in this case, $\pi_n$ can be chosen to
be linear. Unfortunately this assumption does not cover
$\ell_\infty$ spaces, which satisfy only the weak finite-dimensional
approximation property considered in \cite{ak2}.

We begin with a technical lemma on the commutativity of slice and
push-forward; the validity of this identity for rectifiable currents
is proved in \cite[Lemma~5.9]{ak2}; its extension to $\fc{k}{E}$ can
be proved arguing as in Lemma~\ref{lslice2} and in
Lemma~\ref{lcommu}, so we omit a detailed proof.

\begin{lemma}[Slice and push-forward commute]\label{lcommu2}
Let $f\in {\rm Lip}(E;\R^n)$ and $T\in\fc{k}{E}$.
Let $q:\R^n\to\R^k$, $q^\perp:\R^n\to\R^{n-k}$ be respectively the projections
on the first $k$ coordinates and on the last $n-k$ coordinates.
Then
$$
q^\perp_\sharp\langle f_\sharp T,q,x\rangle=(q^\perp\circ f)_\sharp\langle T,q\circ f,x\rangle
\qquad\text{for $\Leb{k}$-a.e. $x\in\R^k$.}
$$
\end{lemma}

Then, we recall the basic result of \cite{white1}, a consequence of
the deformation theorem in \cite{white2}. Thanks to
Proposition~\ref{pchess} we can state White's result in our language
of currents $\modp$, instead of flat chains with coefficients in
$\Z_p$.

\begin{theorem}\label{tbasicwhite}
Let $T\in\fc{k}{\R^N}$. Then $T=0$ $\modp$ if and only if, for all
orthogonal projections $q$ on a $k$-dimensional subspace of $\R^N$,
$\langle T,q,x\rangle=0$ $\modp$ for $\Leb{k}$-a.e. $x\in\R^k$.
\end{theorem}

\begin{theorem}\label{tabs11}
Assume that $E$ is a compact convex subset of a Banach space $F$
having the strong finite-dimensional approximation property, and
that $E$ is a Lipschitz retract of $F$.\\
Let $T\in\fc{k}{E}$ with finite ${\bf M}_p$ mass and assume that,
for some $m\in [1,k]$ the following property holds:
$$ \text{for all $\pi\in [{\mathcal D}]^m$, $\langle
T,\pi,x\rangle=0$ $\modp$ for $\Leb{m}$-a.e. $x\in\R^m$.}$$ Then
$T=0$ $\modp$.
\end{theorem}

\noindent {\sc Proof.} We shall directly prove the statement in the
case $m=k$, which obviously implies all others, by the definition of
iterated slice operator. Let $\pi_n:F\to F$ be given by the strong
finite-dimensional approximation property and let
$T_n=\pi_{n\sharp}T$. We shall prove in the first step that $T_n=0$
$\modp$, and in the second one that $T_n\to T$ in $\fflat$ distance
in $F$. Considering the images $S_n$ of $T_n$ under a Lipschitz
retraction of $F$ onto $E$, which converge to $T$ in $\fflat$
distance in $E$ and are still equal to $0$ $\modp$, this implies
that $T=0$ $\modp$.

\noindent {\bf Step 1.} Since the range of $\pi_n$ is
finite-dimensional we can obviously think of $T_n$ as a flat chain in
a suitable Euclidean space $\R^N$. So, by Theorem~\ref{tbasicwhite},
it suffices to show that the slices induced by orthogonal
projections $q$ on $k$-planes vanish. With no loss of generality we
can assume that $q$ is the orthogonal projection on the first $k$
coordinates, and apply Lemma~\ref{lcommu2} with $f=\pi_n$ to obtain
that
$$
q^\perp_\sharp \langle T_n,q,x\rangle=0\quad\text{$\modp$}
$$
for $\Leb{k}$-a.e. $x\in\R^k$. But since $\|\langle
T_n,q,x\rangle\|_p$ is concentrated on $\{q=x\}$, and
$q^\perp:\{q=x\}\to\{q=0\}$ is an isometry, it follows that $\langle
T_n,q,x\rangle=0$ $\modp$ for $\Leb{k}$-a.e. $x\in\R^k$.

\noindent {\bf Step 2.} Let $E_1$ be the compact metric space
$E\cup\bigcup_n \pi_n(E_n)$ and let $E_2\subset F$ be its closed
convex hull. In order to conclude, it suffices to show that
$\fflat(T_n-T)\to 0$ in $E_2$. Taking into account subadditivity of
the flat norm and \eqref{amb1}, it suffices to show that
$\fflat(\pi_{n\sharp}R-R)\to 0$ for all $R\in\rc{k}{E_2}$; by
density in mass norm, it suffices to prove this fact for
$R\in\ic{k}{E_2}$. Obviously $\pi_{n\sharp}R\to R$ weakly in $E_2$,
i.e. in the duality with Lipschitz forms; then, it suffices to apply
\cite{Wenger-flatconv} to obtain convergence in flat norm in $E_2$.
\qed

We can now prove two basic absolute continuity properties of
$\|T\|_p$.

\begin{theorem}\label{tabs2}
Let $T\in\fc{k}{E}$ with finite ${\bf M}_p$ mass. Then
$\|T\|_p\ll\|T\|_p^*$. In particular $\|T\|_p$ is concentrated on a
countably $\Haus{k}$-rectifiable set.
\end{theorem}

\noindent {\sc Proof.} We fix a compact set $K$ such that
$\|T\|_p^*(K)=0$ and we have to show that $\|T\|_p(K)=0$. Let
$\pi\in [{\mathcal D}]^k$, $u=d(\cdot,K)$ and let $N\subset\R$ be
the Lebesgue negligible set as in \eqref{fasta}. If we consider any
sequence $(s_j)\subset\R\setminus N$ with $s_j\downarrow 0$, then
$T\res\{u<s_j\}$ converges with respect to $\fflat_p$ to
$S\in\fc{k}{E}$ with $[S]=[T]\res K$. The commutativity of slice and
restriction gives
$$
\|\langle T\res\{u<s\},\pi,x\rangle\|_p(E)=\|\langle
T,\pi,x\rangle\|_p(\{u<s\})\qquad \text{for $\Leb{k}$-a.e.
$x\in\R^k$}
$$
for $\Leb{1}$-a.e. $s>0$. Choosing $(s_j)\downarrow 0$ with this
additional property, and assuming also that
$\sum_j\fflat_p(T\res\{u<s_j\}-S)<\infty$, from \eqref{cachanv} we
infer
$$
\lim_{j\to\infty}\langle T\res\{u<s_j\},\pi,x\rangle=\langle
S,\pi,x\rangle\qquad\text{with respect to $\fflat_p$, for
$\Leb{k}$-a.e. $x\in\R^k$.}
$$
Since
$$
\|\langle T\res\{u<s_j\},\pi,x\rangle\|_p(E)=\|\langle
T,\pi,x\rangle\|_p(\{u<s_j\})\to 0\qquad \text{for $\Leb{k}$-a.e.
$x\in\R^k$}
$$
it follows that $\langle S,\pi,x\rangle=0$ $\modp$ for
$\Leb{k}$-a.e. $x\in\R^k$. Then, Theorem~\ref{tabs11} gives that
$[T]\res K=0$. By \eqref{localK} it follows that
$\|T\|_p(K)=\mass_p([T]\res K)=0$.\qed

A direct consequence of Corollary~\ref{cabs3}, ensuring the absolute
continuity of $\|T\|_p$ with respect to $\|T\|_p^*$, is the density
upper bound
\begin{equation}\label{dlimsup}
\limsup_{r\downarrow 0}\frac{\|T\|_p(B_r(x))}{r^k}<\infty\qquad
\text{for $\|T\|_p$-a.e. $x\in E$.}
\end{equation}

Indeed, general covering arguments imply that the set of points
where the $\limsup$ is $+\infty$ is $\Haus{k}$-negligible (see e.g.
\cite[Theorem~2.4.3]{amtilli}), and hence $\|T\|_p$-negligible.

We are now going to a density lower bound for the measure $\|T\|_p$
that gives, as a byproduct, the finiteness of the measure theoretic
support of flat chains with finite ${\bf M}_p$ mass, completing the
proof of Theorem~\ref{tmain1}. Since we can write $T=R+\partial S$
with $R\in\rc{k}{E}$, possibly replacing $T$ by $T-R$ we need only
to consider chains $T$ with $\partial T=0$.

We use a general principle, maybe first introduced by White
\cite{white}, and then used in \cite{cheeger}, \cite{amb} in
different contexts: any lower semicontinuous and additive energy has
the property that any object with finite energy, when seen on a
sufficiently small scale, is a quasiminimizer.

\begin{proposition}\label{pvitali}
Let $T\in\fc{k}{E}$ with finite $\mass_p$ mass and $\partial
T=0$ $\modp$. Then, for all $\eps>0$ the following holds: for
$\|T\|_p$-a.e. $x$ there exists $r_\eps(x)>0$ such that
\begin{equation}\label{blobbo}
\|T\|_p(\overline{B}_r(x))\leq 2\|S+T\|_p(\overline{B}_r(x))
\end{equation}
 whenever $r\in
(0,r_\eps(x))$, $\|T\|_p(B_r(x))\geq\eps r^k$, $\partial [S]=0$ and
$[S]\res\bigl(E\setminus\overline{B}_r(x)\bigr)=0$.
\end{proposition}

\noindent {\sc Proof.} Assume that for some $\eps>0$ the statement
fails. Then, there exists a compact set $K\subset E$ with
$\|T\|_p(K)>0$ such that, for all $x\in K$, we can find balls
$\overline{B}_r(x)$ with arbitrarily small radius satisfying
$\|T\|_p(\overline{B}_r(x))\geq\eps r^k$ and cycles $[S_{r,x}]$ with
$\|T\|_p(\overline{B}_r(x))\geq 2\|T+S_{r,x}\|_p(\overline{B}_r(x))$
and $[S_{r,x}]\res \bigl(E\setminus\overline{B}_r(x)\bigr)=0$.

Let $\delta>0$. By a classical covering argument (see for instance
\cite[Theorem~2.2.2]{amtilli}), we can find a disjoint family
$\{\overline{B}_{r_i}(x_i)\}_{i\in I}$ of these balls, all with
radii $r_i<\delta$, whose union covers $\Haus{k}$-almost all, and
hence $\|T\|_p$-almost all, of $K$. Set now
$$
[T_\delta]:=[T]+\sum_{i\in I}[S_{r_i,x_i}]=[T]\res \bigl(E\setminus
\cup_{i\in I}\overline{B}_{r_i}(x_i)\bigr)+\sum_{i\in
I}[T+S_i]\res\overline{B}_{r_i}(x_i),
$$
where $S_i:= S_{r_i,x_i}$. By approximating with finite unions and sums, taking into account
that ${\bf M}_p([S_i])\leq 3\|T\|_p(\overline{B}_{r_i}(x_i))/2$, it is
not hard to show that $[T_\delta]\in\fcp{k}{E}$ and that
$$
{\bf M}_p([T_\delta])\leq \|T\|_p(E\setminus\cup_{i\in
I}\overline{B}_{r_i}(x_i)\bigr)+\frac{1}{2}\sum_{i\in
I}\|T\|_p(\overline{B}_{r_i}(x_i))\leq \|T\|_p(E\setminus K)+
\frac{1}{2}\sum_{i\in I}\|T\|_p(\overline{B}_{r_i}(x_i)).
$$
In addition, denoting by $[S_i']$ cycles with $\partial
[S_i']=[S_i]$ given by the isoperimetric inequality we have
$$
\fflat_p([T_\delta]-[T])\leq\fflat_p(\sum_i[S_i])\leq
\fflat_p(\sum_i[S_i'])\leq\delta_k3^{1+1/k}\sum_i\bigl[\|T\|_p(B_{r_i}(x_i))\bigr]^{1+1/k}.
$$
Since $\|T\|_p(B_\delta(x))\to 0$ uniformly on $K$ as
$\delta\downarrow 0$, it follows that $[T_\delta]\to [T]$ in $
\fcp{k}{E}$ as $\delta\downarrow 0$. Hence, letting
$\delta\downarrow 0$ the lower semicontinuity of ${\bf M}_p$
provides the inequality $\|T\|_p(E)\leq\|T\|_p(E\setminus
K)+\|T\|_p(K)/2$, which implies $\|T\|_p(K)=0$.\qed

Then, a general and well-known argument based on the isoperimetric
inequalities and an ODE argument provides the following result:

\begin{theorem}
Let $T\in\fc{k}{E}$ with finite ${\bf M}_p$ mass and $\partial
T=0$ $\modp$. Then
\begin{equation}\label{lbblo}
\liminf_{r\downarrow 0} \frac{\|T\|_p(B_r(x))}{r^k}\geq
c>0\qquad\text{for $\|T\|_p$-a.e. $x$}
\end{equation}
with $c>0$ depending only on $k$. As a consequence $\|T\|_p$ is
concentrated on a countably $\Haus{k}$-rectifiable set with finite
$\Haus{k}$-measure.
\end{theorem}
\noindent {\sc Proof.} First of all, we notice that
$$
\limsup_{r\downarrow 0}\frac{\|T\|_p(B_r(x))}{r^k}>0
\qquad\text{for $\|T\|_p$-a.e. $x\in E$.}
$$
 Indeed, if $B$ is the set where the
$\limsup$ above vanishes, general differentiation results (see e.g.
\cite[Theorem~2.4.3]{amtilli}) imply that $\|T\|_p(B')=0$ for any
Borel set $B'\subset B$ with $\Haus{k}(B')<\infty$. Since, by
Theorem~\ref{tabs2}, $\|T\|_p$ is concentrated on a countably
$\Haus{k}$-rectifiable set, and vanishes on $\Haus{k}$-negligible
sets, it follows that $\|T\|_p(B)=0$.

Let $\eps>0$ and let us denote by $C_\eps$ the set where the
$\limsup$ above is larger than $2\eps$; Proposition~\ref{pvitali}
yields, for $\|T\|_p$-a.e. $x\in C_\eps$, $r_\eps(x)>0$ such that
\eqref{blobbo} holds whenever $r\in (0,r_\eps(x))$,
$\|T\|_p(B_r(x))\geq\eps r^k$ and $[S]$ is a $k$-cycle satisfying
$[S-T]\res \bigl(E\setminus\overline{B}_r(x)\bigr)=0$. Then, the
energy comparison argument and the isoperimetric inequality give,
for all such points $x$, the differential inequality
$$
\frac{d}{dr}\|T\|_p^{1/k}(B_r(x))\geq d_k>0
$$
for $\Leb{1}$-a.e. $r\in (0,r_\eps(x))$ such that
$\|T\|_p(B_r(x))>\eps r^k$, for some constant $d_k$ independent of
$\eps$. Since $x\in C_\eps$ we know that the condition
$\|T\|_p(B_r(x))>\eps r^k$ is fulfilled by arbitrarily small radii
$r\in (0,r_\eps(x))$ and we claim that, provided that $\eps<d_k^k$,
if it holds for some $r$, it holds for all $r'\in (r,r_\eps(x))$:
indeed, if $r'$ is the smallest $r'\in (r,r_\eps(x))$ for which it
fails, in the interval $(r,r')$ the function $\|T\|_p^{1/k}(B_r(x))$
has derivative larger than $d_k$, while $\eps^{1/k}r$ has a smaller
derivative. It follows that $\|T\|_p(B_r(x))>\eps r^k$ for all $r\in
(0,r_\eps(x))$ and the differential inequality yields \eqref{lbblo}
at $\|T\|_p$-a.e. $x\in C_\eps$ with $c=d_k^k$. Since
$\cup_{\eps>0}C_\eps$ cover $\|T\|_p$-almost all of $E$ the proof is
finished.\qed

Finally, we complete the list of announced result with the proof of
Corollary~\ref{cmain3}.

\noindent {\bf Proof.} The statement can be easily checked for
chains $T\in\fc{k}{\R^k}$ since $\fc{k}{\R^k}=\rc{k}{\R^k}$ (recall
that $\R^k$ can't support a nonzero $(k+1)$-dimensional integer
rectifiable current). In the general case, let $T\in\fc{k}{E}$ with
finite ${\bf M}_p$ mass, let $S\subset E$ be a countably
$\Haus{k}$-rectifiable Borel set with finite $\Haus{k}$-measure
where $\|T\|_p$ is concentrated and let $B_i\subset \R^k$ be
compact, $f_i:B_i\to E$ be such that $f_i(B_i)$ are pairwise
disjoint, $\cup_i f_i(B_i)$ covers $\Haus{k}$-almost all of $S$ and
$f_i:B_i\to f_i(B_i)$ is bi-Lipschitz, with Lipschitz constants less
than $2$. By McShane's extension theorem we can also assume that
$f_i:E\to\R^k$ are globally defined and Lipschitz. Then, we can find
$\theta_i\in L^1(\R^k,\Z)$ with $|\theta_i|\leq p/2$ and
$\theta_i=0$ out of $B_i$ such that $(f_i)_\sharp ([T]\res
f_i(B_i))= [\segop\theta_i\segcl]$. Since
$$
\sum_i\int_{B_i}|\theta_i|\,dx= \sum_i{\bf
M}_p(\segop\theta_i\segcl)\leq 2^k\sum_i {\bf M}_p ([T]\res
f_i(B_i))\leq 2^k\|T\|_p(E)<\infty
$$
if follows that the current $S:=\sum_i(f_i)_\sharp
\segop\theta_i\segcl\in\rc{k}{E}$ is well defined and, by
construction, $[S]=[T]$. In addition, since $f_i(B_i)$ are pairwise
disjoint, the multiplicity of $S$ takes values in $[-p/2,p/2]$,
hence \cite[Theorem~8.5]{ambkatz} gives that ${\bf
M}(S)=\mass_p(S)$.\qed

\section{Appendix}\label{s7}

In this appendix we state some technical results.

\begin{lemma}\label{ltechnical1}
Let $k\geq 1$ and $m\in[1,k]$. Let $G:\fc{k-m}{E}\to [0,+\infty]$ be
continuous with respect to the $\fflat$ distance and let
$S\in\rc{k}{E}$. Then, for all $\pi\in\bigl[{\rm Lip}(E)\bigr]^m$,
the function $x\mapsto G(\langle S,\pi,x\rangle)$ is Lebesgue
measurable.
\end{lemma}

\noindent {\sc Proof.} With no loss of generality we can assume that
$\pi\in\bigl[{\rm Lip}_1(E)\bigr]^m$. It suffices to show that
$x\mapsto \langle S,\pi,x\rangle$ is the pointwise $\Leb{m}$-a.e.
limit of simple maps. In order to make a diagonal argument we use,
instead, local convergence in measure, so our goal is to find simple
maps $f_h:\R^m\to\fc{k-m}{E}$ such that
$$
\lim_{h\to\infty}\Leb{m}\left(\{x\in B_R(0):\ \fflat(\langle
S,\pi,x\rangle-f_h(x))>\eps\}\right)=0\qquad\forall
R>0,\,\,\forall\eps>0.
$$
Since $\ic{k}{E}$ is dense in $\rc{k}{E}$, by a first diagonal
argument we can assume with no loss of generality that
$S\in\ic{k}{E}$, so that $S_x:=\langle S,\pi,x\rangle\in\ic{k-m}{E}$ for
$\Leb{m}$-a.e. $x$ and
$$
\int_{\R^m}\mass(S_x)+\mass(\partial
S_x)\,dx\leq\mass(S)+\mass(\partial S)<\infty.
$$

In $\ic{k-m}{E}$ we consider the distance
$$
d(T,T'):=\sup\left\{ |T(fdp)-T'(fdp)|:\ |f|\leq
1,\,\,f,p_1,\ldots,p_{k-m}\in {\rm Lip}_1(E)\right\}.
$$
Notice that $d(T,T')\leq\mass(T-T')$ and $d(\partial T,\partial
T')\leq\mass (T-T')$, so that $d(T,T')\leq\fflat(T-T')$ and $\fflat$
convergence is stronger than $d$-convergence. On the other hand,
thanks to the results in \cite{Wenger-flatconv} the two distances
are equivalent in the sets $\{T\in\ic{k-m}{E}:\
\mass(T)+\mass(\partial T)\leq M\}$, $M>0$. Since, according to
\cite{ak2}, $S_x$ is an $MBV$ map with respect to $d$, we can divide
$\R^m$ in open cubes $Q^j_h$ with sides $1/h$ and apply the
Poincar\'e inequality for $MBV$ maps (see \cite{ambrosiomet} or
\cite{desp}) in each of these cubes to find $x_j\in Q^j_h$ with
$$
\int_{Q^j_h}d(S_{x_j},S_x)\,dx\leq \frac{c}{h}\mu(Q^j_h),
$$
where $\mu$ is the total variation measure of $x\mapsto S_x$. It
follows that the piecewise constant map $g_h$ equal to $S_{x_j}$ on
$Q^j_h$ satisfies $\int_{\R^m} d(g_h(x),S_x)\,dx\rightarrow 0$. In
order to improve this convergence from $d$ to $\fflat$ we argue as
follows: we notice that
$$
\Leb{m}\bigl(\{y\in Q^j_h:\ d(S_y,S_{x_j})>4c
h^{m-1}\mu(Q^j_h)\}\bigr)\leq \frac{1}{4}\Leb{m}(Q^j_h),
$$
\begin{equation}\label{choiceyj}
\Leb{m}\bigl(\{y\in Q^j_h:\ \mass(S_y)+\mass(\partial
S_y)>4h^m\int_{Q^j_h}\mass(S_x)+\mass(\partial S_x)\,dx\}\bigr)\leq
\frac{1}{4}\Leb{m}(Q^j_h)
\end{equation}
and therefore we can find points $y_j$ in the intersection of the
complements of these sets. By the triangle inequality
$$
d(S_{y_j},S_x)\leq 4ch^{m-1}\mu(Q^j_h)+d(S_{x_j},S_x)
$$
so that, if we denote by $f_h$ the piecewise constant map equal to
$S_{y_j}$ on $Q^j_h$, we still have
\begin{equation}\label{brambilla1}
\lim_{h\to\infty}\int_{\R^m}d(f_h(x),S_x)\,dx=0.
\end{equation}
In addition, taking \eqref{choiceyj} into account, we have also
\begin{equation}\label{brambilla2}
\sup_{h\in\N}\int_{\R^m}\mass(f_h(x))+\mass(\partial f_h(x))\,dx\leq
4\int_{\R^m}\mass (S_x)+\mass(\partial S_x)\,dx<\infty.
\end{equation}
By \eqref{brambilla1} and \eqref{brambilla2}, taking into account
that $d$ and $\fflat$ are equivalent on the sets $\{T:\
\mass(T)+\mass(\partial T)\leq M\}$ it is easy to infer the local
convergence in measure of $f_h$ to $S_x$ with respect to $\fflat$
(given $\delta>0$ and $R>0$ it suffices to find $M$ such that all
sets $B_R\cap\{\mass(f_h)+\mass(\partial f_h)>M\}$ and
$B_R\cap\{\mass(f)+\mass(\partial f)>M\}$ and have measure less than
$\delta$, then choose $\eps>0$ such that $d(S,S')<\eps$ implies
$\fflat(S-S')<\delta$ whenever $\mass(S)+\mass(\partial S)\leq M$;
eventually one can use the fact that $B_R\cap\{d(f_h,f)>\eps\}$ has
measure less than $\delta/3$ for $h$ sufficiently large). \qed

We now state a standard result on measurable set-valued functions,
see for instance \cite{castaing}.

\begin{lemma}\label{ltechnical2}
Let us assign for all $x\in\R^k$ a finite set $\Lambda(x)\subset E$,
and let us assume that $\{x:\ \Lambda(x)\cap C\neq\emptyset\}$ is
Lebesgue measurable for all closed sets $C\subset E$. Then the sets
$$
B_n:=\left\{x\in\R^k:\ {\rm card\,}\Lambda(x)=n\right\}
$$
are Lebesgue measurable and there exist Lebesgue measurable maps
$f_{j_1},\ldots,f_{j_n}:B_n\to E$ such that
\begin{equation}\label{reponbn}
\Lambda(x)=\left\{f_{j_1}(x),\ldots,f_{j_n}(x)\right\}
 \qquad\text{for $\Leb{k}$-a.e. $x\in B_n$.}
\end{equation}
\end{lemma}

Finally, we conclude this appendix by comparing $\fflat_p$ with the
``polyhedral'' flat distance $\fflat^P_p$ in \eqref{flatp}.

\begin{proposition} \label{p1maggio} There exists $C=C(n,k)$ satisfying
$$
\fflat^P_p(T)\leq C\fflat_p(T)\qquad\text{for all $T\in\ic{k}{\R^n}$
weakly polyhedral.}
$$
\end{proposition}
\noindent
{\sc Proof.} Denoting in this proof by $c$ a generic constant depending
on dimension and codimension, let us recall the Federer-Fleming
deformation theorem $\modp$: for $\epsilon>0$ given, any $R\in\ic{k}{\R^n}$
can be written as $P+U+\partial Q$, with $P$ polyhedral on the scale
$\epsilon$, $\mass_p(P)\leq c(\mass_p(R)+\epsilon\mass_p(\partial R))$,
$\mass_p(\partial P)\leq c\mass_p(\partial R)$, $\mass(U)\leq c\epsilon\mass_p(\partial R)$
and $\mass_p(Q)\leq c\epsilon\mass_p(R)$. The main observation is that,
in the case when $\partial R$ is weakly polyhedral,
the construction (based on piecewise affine deformations of $R$ on skeleta of lower and
lower dimension, until dimension $k$ is reached) of $P$, $U$ and $Q$ provides us with
a current $U$ which is weakly polyhedral as well. Indeed, $U$ corresponds to the $k$-surface
swapt by $\partial R$ during the deformation.

Now, assume that $T=R+\partial S$ with $R\in\ic{k}{\R^n}$ and $S\in\ic{k+1}{\R^n}$ and let
us write $R=P+U+\partial Q$ as above. Since $\partial R=\partial T$ is weakly polyhedral, it
follows that $U$ is weakly polyhedral as well. Now we write $T=P+U+\partial(S+Q)$ and apply
the deformation theorem again to $S+Q$ to obtain $S+Q=P'+U'+\partial Q'$. Again, since
$\partial (S+Q)=\partial (T-P-U)$ is weakly polyhedral, we know that $U'$ is weakly polyhedral.
Now we have $T=(P+U)+\partial (P'+U')$ where $P+U$ and $P'+U'$ are both weakly polyhedral, so that
$\fflat^P_p(T)\leq\mass_p(P+U)+\mass_p(P'+U')$. We have also
$$
\mass_p(P)\leq c(\mass(R)+\epsilon\mass_p(\partial R))=
c\mass_p(R)+c\epsilon\mass_p(\partial T), \quad \mass_p(U)\leq
c\epsilon\mass_p(\partial R)=c\epsilon\mass_p(\partial T).
$$
Analogously we have
$\mass_p(P')+\mass_p(U')\leq c\mass_p(R)+c\mass_p(S)+c\epsilon(\mass_p(T)+\mass_p(\partial T))$
and, since $\epsilon>0$ is arbitrary, we conclude.
\qed

\end{document}